\documentclass[11pt,a4paper]{amsart}
\usepackage[all]{xy}
\usepackage{amssymb}
\usepackage{amsmath}
\usepackage[colorlinks,linkcolor=blue]{hyperref}
\usepackage{mathrsfs}
\usepackage{setspace} 

\newcommand{\wedgep}{\mathsf{\Lambda}}

\newcommand{\scha}{\mathcal{A}}
\newcommand{\schd}{\mathcal{D}}
\newcommand{\schf}{\mathcal{F}}
\newcommand{\schh}{\mathcal{H}}
\newcommand{\schm}{\mathcal{M}}
\newcommand{\schp}{\mathcal{P}}
\newcommand{\schs}{\mathcal{S}}
\newcommand{\scht}{\mathcal{T}}
\newcommand{\schy}{\mathcal{Y}}

\newcommand{\shb}{\mathscr{B}}
\newcommand{\shf}{\mathscr{F}}
\newcommand{\shg}{\mathscr{G}}
\newcommand{\shh}{\mathscr{H}}
\newcommand{\shi}{\mathscr{I}}
\newcommand{\sho}{\mathscr{O}}
\newcommand{\shu}{\mathscr{U}}
\newcommand{\shv}{\mathscr{V}}

\newcommand{\mukai}{\textup{\textsf{v}}}
\newcommand{\deriver}{\textup{R}\!}
\newcommand{\bbmH}{\textbf{\textup{H}}}
\newcommand{\Ku}{\cate{Ku}}
\newcommand{\DCoh}{\cate{D}^b}
\newcommand{\Num}{\textsf{\textup{N}}}

\DeclareMathOperator{\NS}{NS}
\DeclareMathOperator{\Aut}{Aut}
\DeclareMathOperator{\Hom}{Hom}
\DeclareMathOperator{\rank}{rank}

\DeclareMathOperator{\BN}{BN}
\DeclareMathOperator{\Jac}{Jac}
\newcommand{\CJac}{\overline{\Jac}}
\DeclareMathOperator{\Sing}{Sing}
\DeclareMathOperator{\Stab}{Stab}
\newcommand{\stable}{\textsf{\textup{st}}}
\newcommand{\modulis}{\textsf{\textup{M}}}

\DeclareMathOperator{\Gr}{Gr}
\DeclareMathOperator{\CGr}{CGr}

\newcommand{\bz}{\mathbb{Z}}
\newcommand{\br}{\mathbb{R}}
\newcommand{\bc}{\mathbb{C}}
\newcommand{\bp}{\mathbb{P}}

\newcommand{\ch}{\textup{\textsf{ch}}}
\newcommand{\td}{\textup{\textsf{td}}}
\newcommand{\forg}{\textup{\textsf{forg}}}
\newcommand{\sfinf}{\textup{\textsf{inf}}}
\newcommand{\Mukai}{\tilde{\textup{H}}}

\newcommand{\GL}{\operatorname{GL}}
\newcommand{\GLp}{\operatorname{GL^+}}
\newcommand{\grp}{{\tilde{\GLp}}(2,\br)}

\newcommand{\defi}[1]{{\textbf{\emph{#1}}}}
\newcommand{\cate}[1]{{\textbf{\textup{#1}}}}
\newcommand{\funct}[1]{{\textup{\textbf{\textsf{#1}}}}}

\theoremstyle{plain}
\newtheorem{theorem}{Theorem}[section]

\newtheorem{corollary}[theorem]{Corollary}
\newtheorem{proposition}[theorem]{Proposition}
\newtheorem{conjecture}[theorem]{Conjecture}

\theoremstyle{definition}
\newtheorem{definition}[theorem]{Definition}
\newtheorem{remark}[theorem]{Remark}
\newtheorem{example}[theorem]{Example}

\title[EPW varieties as moduli spaces]{EPW varieties as moduli spaces on ordinary GM surfaces and special GM threefolds}
\author{Ziqi Liu}
\address{Dipartimento di Matematica “F. Enriques”, Università degli Studi di Milano, Via Cesare Saldini 50, 20133 Milano, Italy.}
\email{ziqi.liu@unimi.it}

\author{Shizhuo Zhang}
\address{School of Mathematics, Sun Yat-sen University, Guangzhou 510275, China}
\email{zhangshzh28@mail.sysu.edu.cn}

\subjclass[2020]{Primary 14J28 14J45 18G80}

\keywords{Gushel--Mukai varieties, moduli spaces, derived categories}

\thanks{The first author is a member of the INdAM group GNSAGA (2025-2026). The second author is supported by SYSU Starting Grant No. 34000-12256019.}

\begin{document}

\begin{abstract}
We show that the double dual EPW sextic associated with a strongly smooth Gushel--Mukai surface can be realized as a moduli space of semistable objects on its bounded derived category. Also, we observe that the double dual EPW surface associated with a special Gushel--Mukai threefold can be realized as a moduli space of semistable objects on its Kuznetsov component. Then we discuss extensions of our main results to double EPW sextics and double EPW surfaces and a refinement of a statement of Bayer and Perry about Gushel--Mukai threefolds with equivalent Kuznetsov components, under a mild assumption. 
\end{abstract}

\maketitle

\section{Introduction}
In this article, we realize EPW varieties associated with a strongly smooth Gushel--Mukai surface as moduli spaces of semistable objects on its bounded derived category and the Kuznetsov component of the special Gushel--Mukai threefold over it.

Throughout the article, everything is over the field of complex numbers, categories are essentially small, and functors are additive.

\subsection{GM varieties and EPW varieties}
For $n\geq 2$, a \emph{Gushel--Mukai variety} of dimension $n$ is a smooth intersection $X=\CGr(2,V_5)\cap\bp(W)\cap Q$ where $V_5$ is a five-dimensional vector space, $\CGr(2,V_5)\subset\bp(\wedgep^2V_5\oplus\bc)$ is the cone over $\Gr(2,V_5)$, $W\subset\wedgep^2 V_5\oplus\bc$ is an $(n+5)$-dimensional vector subspace, and $Q\subset\bp(W)$ is a quadratic hypersurface. The ample line bundle $\sho_{\bp(W)}(1)$ induces an ample polarization on $X$ whose top intersection degree is always ten. 

The vector space $V_6:=H^0(\bp(W),\shi_X(2))$ of quadrics in $\bp(W)$ containing a Gushel--Mukai variety $X$ is six-dimensional and one can associate $X$ canonically with a Lagrangian subspace $A$ of the wedge space $\wedgep^3 V_6$ according to \cite{DK18}. Given such a Lagrangian subspace $A$, one can construct an integral normal sextic hypersurface $Y^{\geq 1}_{A}\subset\bp(V_6)$, called the \emph{EPW sextic} associated with the data $(V_6,A)$ as in \cite{O'Grady06,O'Grady13}, and its universal covering $\tilde{Y}^{\geq 1}_{A}$ is called the \emph{double EPW sextic} associated with $(V_6,A)$. According to \cite{O'Grady08}, the annihilator $A^{\perp}\subset \wedgep^3 V^\vee_6$ of the subspace $A$ is also a Lagrangian subspace, so one also has the \emph{double dual EPW sextic} $\tilde{Y}^{\geq 1}_{A^{\perp}}$.

The singular locus of $Y^{\geq 1}_{A}$ is an integral normal Cohen--Macaulay surface $Y^{\geq 2}_{A}$ called the \emph{EPW surface} associated with $(V_6,A)$. There exists a double cover over $Y^{\geq 2}_{A}$ branched on its singular locus due to \cite{DK20c}, which is called the \emph{double EPW surface} and denoted by $\tilde{Y}^{\geq 2}_{A}$. Similarly, one has the associated \emph{dual EPW surface} $Y^{\geq 2}_{A^{\perp}}$ and \emph{double dual EPW surface} $\tilde{Y}^{\geq 2}_{A^{\perp}}$. 

\subsection{Kuznetsov components of threefolds}
It is known (see e.g.~\cite{Ku09}) that one can define for a Gushel--Mukai threefold $X$ the \emph{Kuznetsov component}  
$$\cate{Ku}(X)=\{E\in\cate{D}^b(X)\,|\,\deriver\Hom(\shu_X^\vee,E)=\deriver\Hom(\sho_X,E)=0\}$$
where $\shu_X$ is the pullback of the tautological bundle $\shu_G$ on $\Gr(2,V_5)$.

The numerical Grothendieck group $\Num(\Ku(X))$ is a rank two lattice
$$\Num(\Ku(X))\cong\bz[\kappa_1]\oplus\bz[\kappa_2]$$
such that the two orthogonal classes $\kappa_1$ and $\kappa_2$ satisfy
$$\ch(\kappa_1)=-1+\frac{1}{5}h^2\quad\textup{and}\quad \ch(\kappa_2)=2-h+\frac{1}{12}h^3$$
where $h$ is the class $c_1(\sho_X(1))$. One notices that $\chi(\kappa_i,\kappa_i)=-1$ for $i=1,2$.

\subsection{The equivariant structure}
A Gushel--Mukai surface $S=\CGr(2,V_5)\cap\bp(W)\cap Q$ is said to be \emph{strongly smooth} if the intersection $M_S:=\CGr(2,V_5)\cap\bp(W)$ is smooth. In this case, one has $M_S=\Gr(2,V_5)\cap\bp(W)$ and the double cover $X\rightarrow M_S$ branched on $S$ is a Gushel--Mukai threefold, called the \emph{special} Gushel--Mukai threefold over $S$. 

By virtue of \cite{DK18}, the Lagrangian subspace $A(X)$ associated with $X$ is canonically isomorphic to $A(S)$. Hence in this article the notations $A(X)$ and $A(S)$ will be used interchangeably.

The covering involution on $X$ induces an $\mathfrak{S}_2$-group action on $\cate{Ku}(X)$ such that the equivariant category is equivalent to $\cate{D}^b(S)$ according to \cite[Section 8.2]{KP17}. The dual cyclic group $\mathfrak{S}^\vee_2$ acts on $\DCoh(S)$ by an involutive autoequivalence $\Pi$ on $\cate{D}^b(S)$ and the associated equivariant category is canonically equivalent to $\Ku(X)$ by \cite[Theorem 1.3]{Ela15}. 

According to \cite[Section 6.1]{Liu26}, one can find a stability condition $\sigma_S$ on $\DCoh(S)$ for any strongly smooth Gushel--Mukai surface $S$ such that the moduli spaces of semistable objects with respect to $\sigma_S$ detect many geometric features of $S$. Moreover, the stability condition $\sigma_S$ is compatible with the $\mathfrak{S}^\vee_2$-group action on $\cate{D}^b(S)$ so that it induces a stability condition $\sigma_X$ on the equivariant category $\Ku(X)$ and vice versa by virtue of \cite[Section 4.3]{PPZ26}. 

The forgetful functors between $\cate{D}^b(S)$ and $\Ku(X)$ induce covering morphisms between moduli spaces of semistable objects with respect to $\sigma_S$ and $\sigma_X$. The properties and first examples of these covering morphisms have been investigated in \cite{Liu26,PPZ26}. In particular, some classical branched double covers can be reconstructed. This article aims to recover some EPW varieties associated with $S$ and $X$ as moduli spaces, expanding some results in \cite[Section 6]{Liu26}. 

\subsection{Main results}
In \cite{Liu26}, the double sextics $\tilde{Y}^{\geq 1}_{A(S)^{\perp}}$ and $\tilde{Y}^{\geq 1}_{A(S)}$ are realized as moduli spaces of semistable objects with respect to $\sigma_S$ for a Picard rank one Gushel--Mukai surface $S$ based on the descriptions in \cite{O'Grady06,O'Grady15}. This article first generalizes this fact.

The modular realization of $\tilde{Y}^{\geq 1}_{A(S)^{\perp}}$ can be generalized to any strongly smooth Gushel--Mukai surface $S$ using the birational model in \cite[Section 7.5]{DK24} and wall-crossing \cite{BM14-MMP}.

\begin{theorem}\label{00_main-dual-sextic}
The moduli space $\modulis_{\sigma_S}(v_1)$ of semistable objects with Mukai vector $v_1=(1,0,-1)$ is isomorphic to the double dual EPW sextic $\tilde{Y}^{\geq 1}_{A(S)^{\perp}}$ associated with $S$.
\end{theorem}

However, a nice birational model for $\tilde{Y}^{\geq 1}_{A(S)}$ is not generally available. In this case, we have to restrict to the situation that $\tilde{Y}^{\geq 1}_{A(S)}$ is a hyperkähler manifold i.e.~$S$ contains neither lines nor quintic elliptic curves by virtue of \cite{Beri26}. Then we can apply the global Torelli theorem for hyperkähler manifolds to the Hodge isometry induced by the modular realization in the generic case, and use the nef cone description in \cite{BM14-MMP} to conclude the following isomorphism.

\begin{proposition}\label{00_main-sextic-weak}
Suppose that $(S,H)$ is a strongly smooth Gushel--Mukai surface containing no lines or quintic elliptic curves, then the moduli space $\modulis_{\sigma_S}(v_2)$ of semistable objects with Mukai vector $v_2=(-2,H,-2)$ is isomorphic to the double EPW sextic $\tilde{Y}^{\geq 1}_{A(S)}$ associated with $S$.
\end{proposition}

The involution on $\modulis_{\sigma_S}(v_1)$ induced by the $\mathfrak{S}_2^\vee$-action on $\cate{D}^b(S)$ coincides with the canonical involution on the corresponding double dual EPW sextic. Moreover, the ramification locus of $\modulis_{\sigma_S}(v_1)$ is exactly the fixed locus $\modulis_{\sigma_S}(v_1)_{\mathfrak{S}_2^\vee}$ by \cite{DK24,PPZ26}, i.e. $\modulis_{\sigma_S}(v_1)_{\mathfrak{S}_2^\vee}\cong Y^{\geq 2}_{A(X)^{\perp}}$.

According to \cite{Liu26,PPZ26}, the forgetful functor $\Ku(X)\rightarrow\cate{D}^b(S)$ induces a double cover
$$\modulis_{\sigma_X}(\kappa_1)\rightarrow \modulis_{\sigma_S}(v_1)_{\mathfrak{S}_2^\vee}\cong Y^{\geq 2}_{A(X)^{\perp}}$$
branched on the singular locus of the target. 

One can see $\modulis_{\sigma_X}(\kappa_1)\cong\tilde{Y}^{\geq 2}_{A(X)^{\perp}}$ from the description for the moduli space in \cite{JLLZ24} and for the EPW variety in \cite{DK24}. Then one can obtain the statement below following \cite{Liu26}.

\begin{corollary}\label{00_main-dual-surface}
The double cover $\modulis_{\sigma_X}(\kappa_1)\rightarrow \modulis_{\sigma_S}(v_1)_{\mathfrak{S}_2^\vee}$ coincides with $\tilde{Y}^{\geq 2}_{A(X)^{\perp}}\rightarrow Y^{\geq 2}_{A(X)^{\perp}}$.
\end{corollary} 

Similarly, one has a double cover
$$\modulis_{\sigma_X}(\kappa_2)\rightarrow \modulis_{\sigma_S}(v_2)_{\mathfrak{S}_2^\vee}$$
according to \cite{Liu26,PPZ26}, while the parallel result of Corollary~\ref{00_main-dual-surface} for the double EPW surface is not yet available even under the assumption of Proposition~\ref{00_main-sextic-weak}.

In \cite[Theorem 7.7 and Remark 7.9]{DK24}, the authors establish a geometric diagram
$$S^{[2]}\stackrel{g}{\dashrightarrow} X_S\stackrel{f^+}{\longrightarrow}\tilde{Y}_{A(S)^{\perp}}^{\geq 1}$$
where the birational map $g$ is the Mukai flop of all Lagrangian planes in $S^{[2]}$ corresponding to lines in $S$ and the morphism $f^+$ is a symplectic resolution of singularities. In Proposition~\ref{EPW-sextics-dual_the-first-flop} and Proposition~\ref{EPW-sextics-dual_the-second-flop} we recover this diagram by studying the wall-crossing behavior of the moduli space $\modulis_{\sigma_t}(v_1)$ for some stability condition $\sigma_t$ (see Section~\ref{subsection_birational_model}). 

\subsection{Conjecture}
To generalize Proposition~\ref{00_main-sextic-weak} and obtain the analogue of Corollary~\ref{00_main-dual-surface} for double EPW surfaces, we introduce an extra assumption, which we believe is true, so we formulate it as a conjecture. Here one recalls that there exists a canonical polarization $\sho_{\tilde{Y}_{A(X)}^{\geq 1}}(1)$ on $\tilde{Y}_{A(X)}^{\geq 1}$ coming from the primitive ample polarization on the hypersurface $Y_{A(X)}^{\geq 1}$.

\begin{conjecture}\label{00_conjecture}
Consider two Gushel--Mukai threefolds $X$ and $X'$ such that one has an isomorphism $\tilde{Y}_{A(X)}^{\geq 1}\cong\tilde{Y}_{A(X')^{\perp}}^{\geq 1}$ preserving the canonical polarizations, then there exists an equivalence $\cate{Ku}(X)\cong\cate{Ku}(X')$ such that $\kappa_1$ is sent to $\kappa'_2$ and $\kappa_2$ is sent to $\kappa'_1$ up to signs.
\end{conjecture}

This conjecture is a refinement of the homological projective duality \cite[Corollary 6.5]{KP23} for Gushel--Mukai threefolds. It can be proved by assuming a mild relative version of the homological projective duality in \cite{KP21}, which will appear in a subsequent paper. 

Assuming that the conjecture holds, the following results hold. 

\begin{proposition}
Assuming Conjecture \ref{00_conjecture}, then the moduli space $\modulis_{\sigma_S}(v_2)$ of semistable objects is isomorphic to the double EPW sextic $\tilde{Y}^{\geq 1}_{A(S)}$ associated with $S$.
\end{proposition}

\begin{proposition}\label{00_conjecture-generalization2}
Assuming Conjecture \ref{00_conjecture}, then the double cover $\modulis_{\sigma_X}(\kappa_2)\rightarrow \modulis_{\sigma_S}(v_2)_{\mathfrak{S}_2^\vee}$ coincides with $\tilde{Y}^{\geq 2}_{A(X)}\rightarrow Y^{\geq 2}_{A(X)}$.	
\end{proposition}

As a result, one can obtain the following refinement of \cite[Theorem 5.21]{BP23}.

\begin{corollary}\label{00_conjecture-corollary}
Assuming that Conjecture \ref{00_conjecture} holds, and let $X$ and $X'$ be two special Gushel--Mukai threefolds. Then one has a polarized isomorphism $\tilde{Y}^{\geq 1}_{A(X)}\cong \tilde{Y}^{\geq 1}_{A(X')}$ if and only if there is an equivalence $\Ku(X)\cong\Ku(X')$ sending $\kappa_1$ to $\kappa'_1$ and $\kappa_2$ to $\kappa'_2$ up to signs; one has a polarized isomorphism $\tilde{Y}^{\geq 1}_{A(X)}\cong \tilde{Y}^{\geq 1}_{A(X')^{\perp}}$ if and only if there is an equivalence $\Ku(X)\cong\Ku(X')$ sending $\kappa_1$ to $\kappa'_2$ and $\kappa_2$ to $\kappa'_1$ up to signs. 
\end{corollary}

Moreover, one notices that the double EPW surface $\tilde{Y}^{\geq 2}_{A(X)}$ or $\modulis_{\sigma}(\kappa_2)$ is also an invariant for a special Gushel--Mukai threefold $X$. However, it is not clear how to recover the Lagrangian data $(V_6,A)$ from a double EPW surface $\tilde{Y}^{\geq 2}_A$.
 
\section{GM Varieties and EPW Varieties}
In this section, some properties of Gushel--Mukai varieties and the EPW varieties will be recalled. The interested reader can find more information about Gushel--Mukai varieties in the series of articles \cite{DK18,DK19,DK20a,DK20b,DK24} and EPW varieties in \cite{DK20c,O'Grady06,O'Grady08,O'Grady13,O'Grady15,O'Grady16}.

\subsection{The Gushel--Mukai varieties}
In the introduction, we have defined the strongly smooth Gushel--Mukai surfaces and the special Gushel--Mukai threefolds. Here we will provide the general definition. Let $V_5$ be a five-dimensional vector space, and consider the Plücker embedding $\Gr(2,V_5)\subset\bp(\wedgep^2 V_5)$. Let $K$ be a one-dimensional vector space, and consider the cone
$$\CGr(2,V_5)\subset\bp(\wedgep^2V_5\oplus K)$$
with vertex $\nu$. Choose a subspace $W\subset\wedgep^2 V_5\oplus K$ and a quadratic hypersurface $Q\subset\bp(W)$. 

\begin{definition}
A \defi{Gushel--Mukai variety} $X$ is a smooth $n$-dimensional intersection 
$$X=\CGr(2,V_5)\cap\bp(W)\cap Q$$
for $2\leq n\leq 6$, such that the dimension of the vector space $W$ is $n+5$.
\end{definition}

\begin{remark}
This definition (essentially $\dim K=1$) is the same as \cite[Definition 2.1]{DK18} in the smooth case according to \cite[Proposition 2.28]{DK18}.
\end{remark}

A Gushel--Mukai variety $X$ is canonically polarized by restriction of the hyperplane class on $\bp(W)$. The ample polarization $H$ satisfies $H^n=10$ and $K_X=-(n-2)H$. 

\begin{example}\label{GM_surface-condition}
A polarized Gushel--Mukai surface $(S,H)$ is a Brill--Noether general K3 surface of degree $10$. It means that $H$ is very ample, and there are no cubic elliptic curves on $S$.
\end{example}

\begin{example}
The Picard group of a Gushel--Mukai threefold $X$ is generated by the canonical ample polarization $H$, and $H$ makes $X$ a Fano threefold of index $1$ and degree $10$.
\end{example}

The intersection $M_X=\CGr(2,V_5)\cap\bp(W)$ is called the \emph{Grassmannian hull} of the Gushel--Mukai variety $X$. A Gushel--Mukai variety is called \emph{strongly smooth} if its Grassmannian hull is smooth. Any Gushel--Mukai variety of dimension $n\geq 3$ is strongly smooth and the strongly smooth Gushel--Mukai surfaces are characterized by the next statement.

\begin{proposition}[{{\cite[Lemma 2.7, Lemma 2.8]{GLT15}}}]\label{GM_nice-surface-condition}
A Gushel--Mukai surface $(S,H)$ is strongly smooth if and only if it does not have any divisor $E$ satisfying $E^2=0$ and $E\cdot H=4$.
\end{proposition}

In fact, due to \cite[Remark 2.19]{Beri26}, a Gushel--Mukai surface is possibly strongly smooth only when it is ordinary in the following definition.

\begin{definition}
A Gushel--Mukai variety $X$ is called \defi{ordinary} if the vertex of $\CGr(2,V_5)$ is not contained in $M_X$ and is called \defi{special} if the vertex of $\CGr(2,V_5)$ is contained in $M_X$.
\end{definition}

The vertex $\nu$ does not belong to a Gushel--Mukai variety $X$ as $X$ is smooth. So the linear projection from $\nu$ defines a morphism $\gamma\colon X\rightarrow\Gr(2,V_5)$ called the \emph{Gushel map}. It is an embedding when $X$ is ordinary. Suppose that $X$ is special, then it induces a double cover over a smooth intersection $\bp(W/K)\cap \Gr(2,V_5)$, branched along its intersection with a quadric.

It means that an $(n+1)$-dimensional special Gushel--Mukai variety $X$ admits a branched double cover over the Grassmannian hull $M_Y$ of a strongly smooth ordinary Gushel--Mukai variety $Y$ of dimension $n$ such that the branch locus is $Y\subset M_Y$. Conversely, the double cover over $M_Y$ branched on an $n$-dimensional strongly smooth ordinary Gushel--Mukai variety $Y\subset M_Y$ is a special Gushel--Mukai variety of dimension $n+1$.

In this case, the special Gushel--Mukai variety $X$ is called the \emph{opposite} Gushel--Mukai variety of the ordinary Gushel--Mukai variety $Y$ and vice versa.

\subsection{The Lagrangian subspaces and EPW varieties}
Any Gushel--Mukai variety 
$$X=\CGr(2,V_5)\cap\bp(W)\cap Q$$ 
is canonically associated with a Lagrangian subspace $A(X)\subset \wedgep^3 V_6$ for $V_6=H^0(\bp(W),\shi_X(2))$ according to \cite[Theorem 3.10]{DK20a}. One can define a series of projective varieties associated with the Lagrangian subspace. These varieties are introduced and intensively studied by O'Grady.

In general, let $V_6$ be a $6$-dimensional vector space and $A\subset\wedgep^3V_6$ be the Lagrangian subspace for the $\det(V_6)$-valued symplectic form $\omega$ defined by $\omega(\zeta,\eta)=\zeta\land\eta$, then the subset
$$\{[v]\in\bp(V_6)\,|\,\dim(A\cap \left(v\land\wedgep^2V_6\right))\geq\ell\}$$
can be endowed with a scheme structure $Y_A^{\geq \ell}$ for any integral $\ell$ according to \cite[Section 2]{O'Grady06}.

The subscheme $Y_A^{\geq 1}\subset \bp(V_6)$ can a priori be the whole space $\bp(V_6)$. To avoid this, one needs an additional condition. 

\begin{definition}
A non-zero element of a Lagrangian subspace $A\subset\wedgep^3V_6$ is called \defi{decomposable} once it can be written in the form $v_1\land v_2\land v_3$. 
\end{definition}

The Lagrangian subspace associated with a strongly smooth Gushel--Mukai variety does not contain any decomposable element.

\begin{theorem}[{{\cite[Theorem B.2]{DK18}}}]
Consider a Lagrangian subspace $A\subset\wedgep^3V_6$ which does not contain decomposable elements, then
\begin{itemize}
	\item[(a)] $Y_A^{\geq 1}$ is an integral normal sextic hypersurface in $\bp(V_6)$;
	\item[(b)] $Y_A^{\geq 2}=\Sing(Y_A^{\geq 1})$ and is an integral normal Cohen--Macaulay surface of degree $40$; 
	\item[(c)] $Y_A^{\geq 3}=\Sing(Y_A^{\geq 2})$ is finite and smooth;
	\item[(d)] $Y_A^{\geq 4}=\varnothing$.
\end{itemize}
The sextic hypersurface $Y_A^{\geq 1}\subset\bp(V_6)$ is called the EPW sextic associated with $A$.
\end{theorem}

The degree six hypersurface $Y_A^{\geq1}\subset\bp(V_6)$ contains a unique prime ample class $\sho_{Y_A^{\geq1}}(1)$ according to the Lefschetz theorem.

\begin{theorem}[O'Grady \cite{O'Grady06,O'Grady13}]
Consider a Lagrangian subspace $A\subset\wedgep^3V_6$ which does not contain decomposable elements and the associated EPW sextic $Y^{\geq 1}_A\subset\bp(V_6)$, then there exists a double covering $f_A\colon\tilde{Y}^{\geq1}_A\rightarrow Y_A^{\geq 1}$ branched along $Y_A^{\geq 2}$ such that
\begin{itemize}
	\item[(a)] the variety $\tilde{Y}_A^{\geq1}$ is normal and simply connected;
	\item[(b)] the singular locus of $\tilde{Y}_A^{\geq1}$ is the finite set $f^{-1}_A(Y_A^{\geq 3})$;
	\item[(c)] the variety $\tilde{Y}_A^{\geq1}$ is a $K3^{[2]}$-type hyperkähler manifold once $Y_A^{\geq 3}=\varnothing$.
\end{itemize}
The normal variety $\tilde{Y}_A^{\geq 1}$ is called the double EPW sextic associated with $A$. It is canonically polarized by the ample class $\sho_{\tilde{Y}_A^{\geq1}}(1):=f_A^*\sho_{Y_A^{\geq1}}(1)$.
\end{theorem}

One of the most important properties is that the double EPW sextic can be determined by its associated Lagrangian subspace.

\begin{theorem}[O'Grady \cite{O'Grady16}]\label{GM_period-double-EPW}
Consider two Lagrangian subspaces $A\subset\wedgep^3V_6$ and $A'\subset\wedgep^3V'_6$ which contain no decomposable elements, then an isomorphism between $(\tilde{Y}^{\geq 1}_A,\sho_{\tilde{Y}_A^{\geq1}}(1))$ and $(\tilde{Y}^{\geq 1}_{A'},\sho_{\tilde{Y}_{A'}^{\geq1}}(1))$ induces a linear isomorphism $\varphi\colon V_6\rightarrow V'_6$ such that $\wedgep^3\varphi(A)=A'$.
\end{theorem}

The double cover construction has been generalized in \cite{DK20c} to other EPW varieties $Y^{\geq \ell}_A$ and the double EPW surface defined below is particularly interesting.

\begin{theorem}[{{\cite[Theorem 5.2]{DK20c}}}]
Consider a Lagrangian subspace $A\subset\wedgep^3V_6$ which does not contain decomposable elements and the associated EPW surface $Y^{\geq 2}_A\subset\bp(V_6)$, then there exists a unique double covering $f_A\colon\tilde{Y}^{\geq2}_A\rightarrow Y_A^{\geq 2}$ branched over $Y_A^{\geq 3}$ such that
$$f_{A*}\sho_{\tilde{Y}_A^{\geq2}}\cong \sho_{Y_A^{\geq2}}\oplus\omega_{Y_A^{\geq2}}(-3)$$
where $\sho_{Y_A^{\geq2}}(1)$ is induced by $\sho_{\bp(V_6)}(1)$. Moreover, one has
\begin{itemize}
	\item[(a)] the variety $\tilde{Y}_A^{\geq2}$ is integral and normal;
	\item[(b)] the singular locus of $\tilde{Y}_A^{\geq2}$ is the finite set $f^{-1}_A(Y_A^{\geq 3})$;
	\item[(c)] the variety $\tilde{Y}_A^{\geq2}$ is a smooth surface of general type once $Y_A^{\geq 3}=\varnothing$.
\end{itemize}
The variety $\tilde{Y}_A^{\geq 2}$ is called the double EPW surface associated with the Lagrangian subspace $A$.
\end{theorem}

The annihilator $A^{\perp}\subset\wedgep^3V_6^\vee$ of a given Lagrangian subspace $A\subset\wedgep^3V_6$ is also a Lagrangian subspace called the \emph{orthogonal} of $A$. It contains no decomposable elements if and only if $A$ contains no decomposable elements. So one can define the \emph{dual EPW sextic} $Y_{A^{\perp}}^{\geq 1}$, the \emph{double dual EPW sextic} $\tilde{Y}_{A^{\perp}}^{\geq 1}$, the \emph{dual EPW surface} $Y_{A^{\perp}}^{\geq 2}$, and the \emph{double dual EPW surface} $\tilde{Y}_{A^{\perp}}^{\geq 2}$ for a Lagrangian subspace $A\subset\wedgep^3V_6$ as the respective EPW varieties associated with $A^{\perp}$. 

\subsection{The period of Gushel--Mukai threefolds}
Here let us restrict to the case of Gushel--Mukai threefolds. A Gushel--Mukai threefold $(X,H)$ determines a Lagrangian subspace $A$ in $\wedgep^3 V_6$ for the vector space $V_6=H^0(\bp(W),\shi_X(2))$. Together with the vector space $V_5$ in the definition of a Gushel--Mukai variety, one gets a triple $(V_6,V_5,A)$ which in turn determines $X$.

\begin{definition}
A \defi{smooth Lagrangian data set} is a collection $(V_6,V_5,A)$, where
\begin{itemize}
	\item $V_6$ is a six-dimensional vector space;
	\item $V_5\subset V_6$ is a hyperplane;
	\item $A\subset\wedgep^3V_6$ is a Lagrangian subspace with no decomposable elements.
\end{itemize}
An isomorphism of smooth Lagrangian data sets between $(V_6,V_5,A)$ and $(V'_6,V'_5,A')$ is a linear isomorphism
$\varphi\colon V_6\rightarrow V'_6$ such that $\varphi(V_5)=V'_5$ and $(\wedgep^3\varphi)(A)=A'$.
\end{definition}

One has $\dim(A\cap\wedgep^3V_5)=3$ for the smooth Lagrangian data set $(V_6,V_5,A)$ associated with a special Gushel--Mukai threefold, and has $\dim(A\cap\wedgep^3V_5)=2$ for the smooth Lagrangian data set $(V_6,V_5,A)$ associated with an ordinary Gushel--Mukai threefold. Conversely, one can see

\begin{theorem}[Debarre--Kuznetsov \cite{DK18}]
The map sending a Gushel--Mukai threefold to its smooth Lagrangian data induces a bijection between the set of Gushel--Mukai threefolds and the set of smooth Lagrangian data $(V_6,V_5,A)$ such that $\dim (A\cap\wedgep^3V_5)\in\{2,3\}$.
\end{theorem}

Since the double EPW sextic $\tilde{Y}_A^{\geq 1}$ is determined by the data $(V_6,A)$, one cannot distinguish the Gushel--Mukai threefolds by their associated double EPW sextics directly.

\begin{definition}
Two Gushel--Mukai threefolds are called \defi{period partners} when they admit isomorphic polarized double EPW sextics. Two Gushel--Mukai threefolds are called \defi{period duals} when the double EPW sextic of one is the double dual EPW sextic of the other.
\end{definition} 

Moreover, one has the following criterion by \cite[Theorem 3.1 and Lemma 3.8]{KP18}. Here one defines $Y_{A(X)}^{\ell}:=Y_{A(X)}^{\geq\ell}-Y_{A(X)}^{\geq \ell+1}$ for $\ell\geq 1$.

\begin{proposition}\label{GM_threefold-period-dual}
Consider a Gushel--Mukai threefold $X$. Then a point $q$ in $Y_{A(X)}^{\geq 2}$ corresponds to a period dual $X_q$ of $X$. In particular, the period dual
\begin{itemize}
	\item is a special Gushel--Mukai threefold if the point $q$ is in $Y_{A(X)}^3$;
	\item is an ordinary Gushel--Mukai threefold if the point $q$ is in $Y_{A(X)}^2$.
\end{itemize}
Conversely, any period dual of $X$ arises as above for some $q\in\bp(V_6(X))$.
\end{proposition}

Due to \cite[Theorem 1.9]{BP23}, information about the period partners or duals for a Gushel--Mukai threefold $X$ can be read from the Kuznetsov component $\Ku(X)$. Once an additional datum is included (serving as $V_5$), one can recover $X$ from $\Ku(X)$ according to \cite{JLZ26}.

\section{Stability conditions, moduli spaces, and on K3 surfaces}
\subsection{Stability conditions on triangulated categories}
To reconstruct EPW varieties, we need to use the notion of stability conditions introduced by Bridgeland \cite{Bri07}. Here we will adopt a more recent formulation of stability conditions proposed in \cite[Appendix 1]{BMS16} based on the locally-finite full stability conditions in \cite{Bri07,Bri08}.

Here one fixes a triangulated category $\schd$, and a group homomorphism $\mukai\colon K(\schd)\rightarrow\Lambda$ from the Grothendieck group of $\schd$ to a finite rank lattice $\Lambda$.

\begin{definition}
A \defi{stability condition} on $\schd$ is a pair $\sigma=(\schp,Z)$ where
\begin{itemize}
	\item the \emph{central charge} $Z$ is a linear map $Z\colon\Lambda\rightarrow\bc$, and
	\item the \emph{slicing} $\schp$ is a collection of full subcategories $\schp(\phi)\subset\schd$ for all $\phi\in\br$,
\end{itemize}
such that
\begin{itemize}
	\item $\schp(\phi+1)=\schp(\phi)[1]$ and $\Hom(\schp(\phi_1),\schp(\phi_2))=0$ for any $\phi_1>\phi_2$;
	\item one has $Z(\mukai(E))\in\br_{>0}\cdot\exp(\pi\phi\sqrt{-1})$ for any $0\neq E\in\schp(\phi)$;
	\item every object $E$ of $\schd$ admits a Harder--Narasimhan filtration
	\begin{displaymath}
			0= \xymatrix{
				E_0\ar[r]&E_1\ar[r]\ar[d]&E_2\ar[r]\ar[d]&\cdots\ar[r]&E_{n-1}\ar[r]\ar[d]&E_n\ar[d]\\
				&A_1\ar@{-->}[ul]&A_2\ar@{-->}[ul]&\ar@{-->}[ul]\cdots&\cdots&A_n\ar@{-->}[ul]
			} =E
	\end{displaymath}
	with $A_i\in\schp(\phi_i)$ for all $i$ such that $\phi_1>\cdots>\phi_n$. 
	\item there exists a quadratic form $Q$ on the vector space $\Lambda_{\br}$ such that $\ker(Z_{\br})$ is negative definite with respect to $Q$, and one has $Q(\mukai(E))\geq0$ for any $E\in\schp(\phi)$ and any $\phi\in\br$.
\end{itemize}
In this case, each $\schp(\phi)$ is an abelian category. A non-zero object of $\schd$ is called \defi{semistable} with phase $\phi$ once it belongs to $\schp(\phi)$ and is called \defi{stable} once it is semistable and simple in $\schp(\phi)$.
\end{definition}

\begin{remark}
Let $\sigma$ be a stability condition, then the Harder--Narasimhan filtration for an object $E$ is unique up to isomorphism. In this case, the objects $A_i$ are called the \emph{factors} of $E$. 
\end{remark}

\begin{remark}
Let $\sigma=(\schp,Z)$ be a stability condition, then any abelian category $\schp(\phi)$ has finite length (see for example \cite[page 330]{Bri07}) and an object in $\schp(\phi)$ therefore admits a Jordan--Hölder filtration with uniquely determined Jordan--Hölder factors. 
\end{remark}

\begin{remark}
Let $\sigma=(\schp,Z)$ be a stability condition on $\schd$, then the smallest extension-closed subcategory $\schp((0,1])\subset\schd$ containing all the abelian subcategories $\schp(\phi)$ for $\phi\in(0,1]$ is the heart of a bounded $t$-structure on $\schd$. In this case, it is called the \emph{heart} of $\sigma$.
\end{remark}

This article will consider only certain stability conditions on triangulated categories arising from geometry.

\begin{definition}
A triangulated category $\schd$ is called \defi{geometric} once it is equivalent to a triangulated subcategory of the bounded derived category of a smooth projective variety $X$ such that the embedding $\schd\hookrightarrow\DCoh(X)$ has both the left and the right adjoint functors.
\end{definition}

The Euler form $\chi$ on a smooth projective variety $X$ restricts to a geometric triangulated category $\schd\subset\DCoh(X)$ and one defines $\Num(\schd):=K(\schd)/\ker(\chi)$.

\begin{definition}
A stability condition on a geometric triangulated category is called \defi{numerical} once the surjection $\mukai$ is the quotient $K(\schd)\twoheadrightarrow \Num(\schd):=K(\schd)/\ker(\chi)$. 
\end{definition}

Henceforth, a triangulated category $\schd$ is always geometric and a stability condition is always numerical. According to \cite{Bri07}, the set of (numerical) stability conditions on $\schd$, if non-empty, admits a complex manifold structure $\Stab(\schd)$ with dimension equal to the rank of $\Num(\schd)$.

Also, there are two mutually-commutative group actions on the complex manifold $\Stab(\schd)$ due to \cite[Lemma 8.2]{Bri07} as follows.

The group $\Aut(\schd)$ acts on $\Stab(\schd)$ on the left, via the mapping
$$\Phi. (\schp, Z):= (\schp', Z \circ \Phi_*^{-1})  $$
where $\Phi_*$ is the induced automorphism on $\Num(\schd)$ and $\schp'$ is defined by $\schp'(\phi):=\Phi(\schp(\phi))$.

The universal cover 
$$
\grp  = \left \{ \tilde{g}=(M, f) ~  \middle\vert ~ 
\begin{array}{l} 
	M \in \GLp(2,\br), \text{ and }f \colon \br \rightarrow \br \text{ is an increasing} \\  
	\text{function such that for all $\phi \in \br$ we have }f(\phi+1)  \\
	=f(\phi) + 1 \text{ and $M \cdot e^{i\pi \phi} \in \br_{>0} \cdot e^{i \pi f(\phi)}$} 
\end{array}  \right \}
$$
of the group $\GLp(2,\br)=\{M\in\GL(2,\br)\,|\,\det(M)>0\}$ acts on the right on $\Stab(\schd)$, via 
$$
(\schp=\{\schp(\phi)\}_{\phi\in\br}, Z) . (M, f):= (\schp'=\{\schp(f(\phi))\}_{\phi\in\br}, M^{-1} \circ Z). 
$$ 
where we identify $\bc$ with $\br\oplus\br\sqrt{-1}$ to make sense of the composite $M^{-1}\circ Z\colon \Num(\schd)\rightarrow\bc$.

\subsection{Moduli spaces of semistable objects}
Let $\sigma$ be a stability condition on a triangulated category $\schd\subset\DCoh(X)$ and $v$ be an element in $\Num(\schd)$, then one can define a moduli stack $\mathfrak{M}_{\sigma}(v)$ of $\sigma$-semistable objects with numerical class $v$. A comprehensive reference is \cite{BLMNPS21}.

\begin{definition}
The moduli stack $\mathfrak{M}_{\sigma}(v)\colon(\cate{Sch}/\bc)^{\cate{op}}\rightarrow\cate{Gpds}$ is defined by setting 
$$
\mathfrak{M}_{\sigma}(v)(T)= \left \{ E\in\cate{D}^b(X_T) ~  \middle\vert ~ 
\begin{array}{l} 
	\,E|_{X_t}\textup{ is in }\schd_t\subset\cate{D}^b(X_t)\textup{ and} \\ 
	\textup{is }\sigma\textup{-semistable with class }v \\
	\textup{for any closed point }t\in T 
\end{array}  \right \}
$$
for any locally finitely generated complex scheme $T$. Similarly, one defines the moduli stack $\mathfrak{M}^{\stable}_{\sigma}(v)$ of stable objects, which is a substack of $\mathfrak{M}_{\sigma}(v)$.
\end{definition}

A stability condition is called \emph{proper} once, for any $v$, the moduli stack $\mathfrak{M}_{\sigma}(v)$ is an algebraic stack of finite type, and there is an open immersion $\mathfrak{M}_{\sigma}^{\stable}(v)\subset \mathfrak{M}_{\sigma}(v)$, and $\mathfrak{M}_{\sigma}(v)$ admits a proper good moduli space $\modulis_{\sigma}(v)$ in the sense of \cite{Al13}, and the open substack $\mathfrak{M}_{\sigma}^{\stable}(v)$ is a $\mathbb{G}_m$-gerbe over the counterpart $\modulis^{\stable}_{\sigma}(v)\subset\modulis_{\sigma}(v)$. All the stability conditions in this article are proper.

Such a good moduli space $\modulis_{\sigma}(v)$ is an algebraic space and the (closed) points in $\modulis_{\sigma}(v)$ are $S$-equivalence classes of semistable objects with phase ranging in $(-1,1]$.

Moreover, the moduli space $\modulis_{\sigma}(v)$ varies with respect to the variation of the stability condition $\sigma$ in the complex manifold $\Stab(\schd)$.

\begin{proposition}[{{\cite[Section 9]{Bri08} and \cite[Proposition 3.3]{BM11}}}]\label{K3-Stab_definition-of-wall}
Given a primitive $v\in \Num(\schd)$, then there exists a locally finite set of walls (real codimension one submanifolds with boundary) in the complex manifold $\Stab(\schd)$, depending only on $v$, such that:
\begin{itemize}
	\item The sets of $\sigma$-stable objects of class $v$ remain unchanged for any stability condition $\sigma$ in the same connected chamber;
	\item When $\sigma$ lies on a single wall in $\Stab(\schd)$, then there is a $\sigma$-semistable object that is unstable in one of the adjacent chambers, and semistable in the other adjacent chamber.
\end{itemize}
In particular, a stability condition $\sigma$ lies on a wall with respect to $v$ only if there exists a strictly $\sigma$-semistable object of class $v$. 
\end{proposition}

In general, the classification of walls is not available but one always has the following two situations for a stability condition.

\begin{definition}
Let $v$ be a class in $\Num(\schd)$, a stability condition $\sigma$ on a triangulated category $\schd$ is called \defi{$v$-generic} when it is not contained in any wall with respect to $v$.
\end{definition}

\begin{definition}
Let $W\subset\Stab(\schd)$ be a wall with respect to a given class $v\in \Num(\schd)$ such that one has $\modulis^{\stable}_{\sigma}(v)\neq\varnothing$ for some $\sigma\in\Stab(\schd)$ close to $W$, then $W$ is called a \defi{totally semistable wall} if one has $\modulis^{\stable}_{\sigma_0}(v)=\varnothing$ for any $\sigma_0\in W$ which is not in any other wall.
\end{definition}

The walls and wall-crossing behaviors of the moduli space $\modulis_{\sigma}(v)$ are well-understood for some stability conditions on $\bp^2$ \cite{LZ19}, K3 surfaces \cite{BM14,BM14-MMP,MZ16} and Enriques surfaces \cite{NY20}. Since a strongly smooth Gushel--Mukai surface is in particular a K3 surface, the article uses some results regarding wall-crossing for K3 surfaces which are recollected in what follows.

\subsection{Some stability conditions on K3 surfaces}
A connected component of the complex manifold $\Stab(\DCoh(S))$ for a K3 surface $S$ is introduced in \cite{Bri08}. In this case, the Mukai vector
$$v(E)=\ch(E).\sqrt{\td(S)}=(\rank(E),c_1(E),\rank(E)+c_1(E)^2/2-c_2(E))$$
identifies the numerical Grothendieck group $\Num(S):=\Num(\DCoh(S))$ with the algebraic Mukai lattice
$$\Mukai(S,\bz)=H^0(S,\bz)\oplus\NS(S)\oplus H^4(S,\bz)$$ 
where $\NS(S)$ is the Néron--Severi group of the K3 surface. The lattice $\Mukai(S,\bz)$ is endowed with the Mukai pairing $\langle(r_1,\Delta_1,s_1),(r_2,\Delta_2,s_2)\rangle=\Delta_1.\Delta_2-r_1.s_2-r_2.s_1$ where $\Delta_1\cdot\Delta_2$ is the usual intersection product in $H^2(S,\bz)$ and carries a canonical weight two Hodge structure from the one on the integer cohomology group $H^0(S,\bz)\oplus H^2(S,\bz) \oplus H^4(S,\bz)$. 

An ample class $\omega\in \NS(S)_{\br}$ defines a notion of $\mu_{\omega}$\emph{-slope stability} on $\cate{Coh}(S)$, by declaring the slope $\mu_{\omega}(\shf):=c_1(\shf)\cdot\omega/\rank(\shf)$ for a torsion-free sheaf $\shf$ and $\mu_{\omega}(\shf)=+\infty$ for a torsion $\shf$. Truncating the HN filtrations at the number $\beta\cdot\omega$, a pair $(\omega,\beta)\in\NS(S)^2_{\br}$ with $\omega$ ample determines a torsion pair $(\scht_{\omega,\beta},\schf_{\omega,\beta})$ on $\cate{Coh}(S)$ such that
\begin{align*}
	\scht_{\omega,\beta}&=\{\textup{torsion-free part of }\shf\textup{ has }\mu_{\omega}\textup{-semistable HN factors of slope} >\beta\cdot\omega\}\\
	\schf_{\omega,\beta}&=\{\shf \textup{ is torsion free and its }\mu_{\omega}\textup{-semistable HN factors have slope} \leq\beta\cdot\omega\}
\end{align*}
This torsion pair gives a bounded $t$-structure on $\cate{D}^b(S)$ with heart 
$$\scha_{\omega,\beta}=\{E\in \cate{D}^b(S)\,|\, \shh^i(E)=0\textup{ for }i\notin\{-1,0\},\shh^{-1}(E)\in\schf_{\omega,\beta},\shh^0(E)\in\scht_{\omega,\beta}\}$$
according to \cite{HRS}, on which one can define a stability condition $\sigma_{\omega,\beta}$ under certain conditions.

\begin{definition}
An object $E$ in $\cate{D}^b(S)$ is called \defi{spherical} once $\deriver\Hom(E,E)=\bc[0]\oplus\bc[-2]$.
\end{definition}

\begin{proposition}[{{\cite[Lemma 6.2 and Section 7]{Bri08}}}]
Consider a K3 surface $S$ and two classes $\omega,\beta\in\NS(S)_{\br}$ with $\omega$ ample, then there exists a stability condition $\sigma_{\omega,\beta}$ with central charge
$$Z_{\omega,\beta}\colon \Mukai(S,\bz)\rightarrow\bc,\quad (r,\Delta,s)\mapsto \frac{1}{2}(2\beta\cdot\Delta-2s+r(\omega^2-\beta^2))+(\Delta-r\beta)\cdot\omega\sqrt{-1}$$
and heart $\scha_{\omega,\beta}$ if and only if $Z_{\omega,\beta}(\shf)\notin\br_{\leq0}$ for any spherical sheaf $\shf$ on $S$.
\end{proposition}

\begin{example}\label{K3-Stab_the-stability-condition-on-GM-surface}
Let $(S,H)$ be a strongly smooth Gushel--Mukai surface endowed with the canonical degree $10$ polarization, then one can check that $(\omega,\beta)=(1/5 H,-2/5 H)$ defines a stability condition $\sigma_S:=\sigma_{\omega,\beta}$ on $\DCoh(S)$ with central charge
$$Z_{\sigma_S}(r,\Delta,s)=-\left(\frac{3}{5}r+\frac{2}{5}H\cdot\Delta+s\right)+\left(\frac{4}{5}r+\frac{1}{5}\Delta\cdot H\right)\sqrt{-1}$$
using the above criterion and Proposition~\ref{GM_nice-surface-condition}. Otherwise, one has $Z_{\omega,\beta}(w)\leq0$ for some spherical class $w=(r,\Delta,s)$ with $r>0$. In particular, one has
$$r>0,\quad\Delta^2=2rs-2,\quad \Delta\cdot H=-4r,\quad 3r+2\Delta\cdot H+5s\geq0$$
at the same time. Set $x=-r$ and $y=\Delta^2/2$, then the only possible integer solutions for $(x,y)$ are $(-1,0)$ and $(-2,3)$. It follows that one has either $H\cdot \Delta=4$ and $\Delta^2=0$ or $H\cdot \Delta=2$ and $\Delta^2=0$. Neither is possible due to Proposition~\ref{GM_nice-surface-condition}.
\end{example}

The stability conditions $\sigma_{\omega,\beta}$ on a K3 surface $S$ are proper and contained in a connected component $\Stab^{\dag}(S)$ of $\Stab(\DCoh(S))$ according to \cite{Bri08,To08}. Moreover, the following statement asserts that, for a certain vector $v\in \Mukai(S,\bz)$ and a given polarization $H$ on $S$, the moduli space $\modulis_{\sigma}(v)$ is equal to the moduli space $\modulis_H(v)$ of Gieseker semistable sheaves on $(S,H)$ with Mukai vector $v$ once the stability condition $\sigma$ belongs to a special chamber in $\Stab^{\dag}(S)$.

\begin{proposition}[\cite{BM14,Bri08,To08}]\label{K3-Stab_K3-large-volume-limit}
Consider a polarized K3 surface $(S,H)$ and a stability condition $\sigma_{
\omega,\beta}$ on $\cate{D}^b(S)$ such that $\omega,\beta\in\br[H]$, then for any given $v\in \Mukai(S,\bz)$ with $\modulis_H(v)\neq\varnothing$ there exists a real number $t_0=t(v)$ such that $\modulis_{\sigma_{t\omega,\beta}}(v)=\modulis_H(v)$ for any $t\geq t_0$.
\end{proposition}

In general, the possible walls in $\Stab^{\dag}(S)$ with respect to a given primitive Mukai vector are classified in \cite{BM14} and the criteria to determine a wall can be found in \cite{BM14-MMP}. These results will be recalled in the subsequent subsection.

\subsection{The walls and wall-crossing on K3 surfaces}\label{K3-Stab_BM-criterion}
Let $S$ be a K3 surface and $v\in\Mukai(S,\bz)$ be a primitive Mukai vector with $\langle v,v\rangle\geq -2$. Suppose that $W\subset\Stab^{\dag}(S)$ is a wall with respect to $v$ in the sense of Proposition~\ref{K3-Stab_definition-of-wall} and $\sigma_0$ is a generic point on $W$, then for any two $v$-generic stability conditions $\sigma_{\pm}$ in the adjacent chambers with respect to $W$ one has a birational map
$$\modulis_{\sigma_-}(v)\dashrightarrow \modulis_{\sigma_+}(v)$$
and two contractions $\modulis_{\sigma_{\pm}}(v)\rightarrow\bar{M}_{\pm}$ due to \cite[Theorem 1.4]{BM14}. 

To determine a potential wall $W$, one introduces a rank two primitive sublattice $\bbmH_W$ of the Mukai lattice $\Mukai(S,\bz)$ such that $w\in\bbmH_W$ if and only if $\Im(Z(w)/Z(v))=0$ for any $\sigma\in W$.

\begin{proposition}[{{\cite[Proposition 5.1]{BM14-MMP}}}]
In the setups of this subsection, the Mukai vector of any HN factor of an object $A$ in $\modulis_{\sigma_+}(v)$ with respect to $\sigma_-$ is contained in $\bbmH_W$, and vice versa. The same holds for any object in $\modulis_{\sigma_0}(v)$. Similarly, any Jordan--Hölder factor of an object $A\in\modulis^{\stable}_{\sigma_0}(w)$ for some $w\in \bbmH_W$ will have Mukai vector contained in $\bbmH_W$.
\end{proposition}

In this case, using the existence of certain semistable objects, one can determine the type of the wall $W$ according to specific vectors in $\bbmH_W$ by \cite[Theorem 5.7]{BM14-MMP} when $\langle v,v\rangle>0$.

\begin{definition}
Let $\langle v,v\rangle> 0$, then the wall $W$ is called a \defi{fake wall} once there are no curves in $\modulis_{\sigma_{\pm}}(v)$ that are $S$-equivalent to each other with respect to $\sigma_0$; a \defi{flopping wall} once $\bar{M}_{+}=\bar{M}_{-}$ and the induced birational map $\modulis_{\sigma_-}(v)\dashrightarrow \modulis_{\sigma_+}(v)$ induces a flopping contraction; a \defi{divisorial wall} once the induced contractions $\modulis_{\sigma_{\pm}}(v)\rightarrow\bar{M}_{\pm}$ are divisorial contractions.
\end{definition}

\begin{theorem}
In the setups of this subsection, and assume that $\langle v,v\rangle> 0$, then the potential wall $W$ is a totally semistable wall if and only if there exists either a vector $w\in\bbmH_W$ such that $\langle w,w\rangle=0$ and $\langle v,w\rangle=1$, or such that $\langle w,w\rangle=-2$, $\Re(Z(w)/Z(v))>0$ and $\langle v,w\rangle<0$. 

In addition, $W$ is a wall inducing a divisorial contraction if and only if there exists a vector $w$ in $\bbmH_W$ such that either $\langle w,w\rangle=0$ and $\langle v,w\rangle\in\{1,2\}$ or $\langle w,w\rangle=-2$ and $\langle v,w\rangle=0$.
\end{theorem}

Otherwise $W$ is either a flopping wall with respect to $v$ or not a wall. In this article, only vectors $v$ with $\langle v,v\rangle=2$ will be involved, so one has the following special criterion.

\begin{theorem}\label{K3-Stab_wall-criteria-flop}
Suppose that $W$ is not a totally semistable wall or a divisorial wall and the vector $v$ satisfies $\langle v,v\rangle=2$, then $W$ is a flopping wall if and only if there exists a vector $w\in\bbmH_W$ such that $\langle w,w\rangle=-2$ and $\langle v,w\rangle=1$.
\end{theorem}

In this case, according to \cite[Proposition 9.4 and Section 14]{BM14-MMP}, the Jordan--Hölder factors for a strictly semistable object $A$ in $\modulis_{\sigma_0}(v)$ are spherical $\sigma_{\pm}$-stable objects. Moreover, one notices that the vector $w'=v-w$ also satisfies $\langle w',w'\rangle=-2$ and $\langle v,w'\rangle=1$ in the setting of Theorem~\ref{K3-Stab_wall-criteria-flop}, so a general flopping wall is induced by two spherical stable objects.

To identify the spherical factors for such a flopping wall, one also needs to understand the walls for $\modulis_{\sigma_0}(v)$ with respect to spherical classes.

\begin{proposition}[{{\cite[Proposition 4.8]{Bo24}}}]\label{K3-Stab_wall-criteria-spherical}
Consider a spherical class $v$ and $A\in\modulis_{\sigma_+}(v)$ which is not $\sigma_0$-stable, then the set of Jordan--Hölder factors for $A$ contains only two spherical objects.
\end{proposition}

\section{The passage to special Gushel--Mukai threefolds}
\subsection{Equivariant categories}
Given an action of $G$ on $\schd$, one can construct the \emph{equivariant category} $\schd_G$ with respect to the $G$-action on $\schd$ whose objects and morphisms are respectively the $G$-equivariant objects and $G$-equivariant morphisms according to \cite{Ela15,BO23}. 

Suppose in addition that $G$ is abelian, then there exists a natural action of the dual group $G^\vee:=\Hom(G,\bc^\times)$ on $\schd_G$ induced by a $G$-action on $\schd$. 

\begin{theorem}[{{\cite[Theorem 4.2]{Ela15}}}]\label{En_Elagin-reconstruction}
Consider a finite abelian group $G$ acting on an idempotent complete category $\schd$, then there exists an equivalence $\funct{X}\colon(\schd_G)_{G^\vee}\cong\schd$ such that the $G$-action is identified with the $(G^\vee)^\vee$-action by $(G^\vee)^\vee\cong G$.
\end{theorem}

There are mutually left and right adjoint functors: the \emph{forgetful functor} 
$$\funct{Forg}\colon\schd_G\rightarrow\schd$$ 
defined by forgetting the $G$-equivariant structure, and the \emph{inflation functor} 
$$\funct{Inf}\colon\schd\rightarrow\schd_G$$ 
given on objects by sending $A$ to $\bigoplus_{g \in G} g_*A$ with the natural $G$-equivariant structure (one can see, for example, \cite[Lemma 3.8]{Ela15}). Moreover, one has the following observation:  

\begin{proposition}\label{En_Elagin-reconstruction-ppst}
Under the assumptions of Theorem~\ref{En_Elagin-reconstruction}, there is an equivalence of functors
$$
\funct{Forg} \cong \funct{Inf} \circ \funct{X} \colon (\schd_G)_{G^\vee} \rightarrow \schd_G
$$
where $\funct{Forg}$ is the one for the induced $G^\vee$-action on $\schd_G$.
\end{proposition}

Suppose in addition that both $\schd$ and $\schd_G$ are triangulated categories, then $(\schd_G)_{G^\vee}$ comes with a natural class of exact triangles, consisting of those triangles whose images under 
$$\funct{Forg}\colon (\schd_G)_{G^\vee} \rightarrow \schd_G$$ 
are exact. Then, according to Proposition~\ref{En_Elagin-reconstruction-ppst}, the equivalence $\funct{X}$ in Theorem~\ref{En_Elagin-reconstruction} respects exact triangles and is therefore an equivalence of triangulated categories.

\subsection{Equivariant stability conditions and covering maps}
Let $G$ be a finite group action on a triangulated category $\schd$. Then one can construct new stability conditions on $\schd_G$ out of the stability conditions on $\schd$ under suitable compatible conditions. To avoid unnecessary technical notions, we will only consider the case where $\schd=\DCoh(X)$ for a smooth projective variety $X$ and the group $G$ is the cyclic group $\mathfrak{S}_2$. One consults \cite{PPZ26} for more general statements.

\begin{theorem}[\cite{MMS09,PPZ26,Pol07}]\label{En_equivariant-stability-condition}
Suppose that the involution $\Pi$ on $\DCoh(X)$ is a generator of the group action $\mathfrak{S}_2$ on $\schd=\DCoh(X)$, and let $\sigma=(\schp,Z)$ be a stability condition on $\schd$ such that for each $\phi\in\br$ one has $\Pi(\schp(\phi))=\schp(\phi)$ and $Z\circ\Pi_*=Z$ where $\Pi_*$ is the action of $\Pi$ on $\Num(X)$, then 
	\begin{align*}
		\schp_{\mathfrak{S}_2}(\phi)&:=\{A\in\schd_{\mathfrak{S}_2}\,|\,\funct{Forg}(A)\in\schp(\phi)\}\\
		Z_{\mathfrak{S}_2}&=Z\circ\forg\colon \Num(\schd_{\mathfrak{S}_2})\rightarrow \bc
	\end{align*}
determines a stability condition $\sigma_{\mathfrak{S}_2}$ on $\schd_{\mathfrak{S}_2}$.
\end{theorem}

In the construction from Theorem~\ref{En_equivariant-stability-condition}, both the central charge $Z_{\mathfrak{S}_2}$ and the slicing $\schp_{\mathfrak{S}_2}$ of the stability condition $\sigma_{\mathfrak{S}_2}$ are preserved by the dual group $\mathfrak{S}_2^\vee$ action on $\schd_{\mathfrak{S}_2}$. So we may apply
the theorem again to obtain a stability condition $(\sigma_{\mathfrak{S}_2})_{\mathfrak{S}_2^\vee}$ on $(\schd_{\mathfrak{S}_2})_{\mathfrak{S}_2^\vee}$. This iterated construction is compatible with the equivalence in Theorem~\ref{En_Elagin-reconstruction}, up to a multiplication by $2$.

\begin{proposition}
In the setup of Theorem~\ref{En_equivariant-stability-condition}, the stability condition $(\sigma_{\mathfrak{S}_2})_{\mathfrak{S}_2^\vee}$ corresponds to the stability condition $2\cdot\sigma=(\schp,2Z)$ under the equivalence $\funct{X}\colon (\schd_{\mathfrak{S}_2})_{\mathfrak{S}_2^\vee}\cong\schd$ of Theorem~\ref{En_Elagin-reconstruction}.
\end{proposition}

\subsection{The Kuznetsov components and the involutions}
Let $X$ be a Gushel--Mukai threefold endowed with the canonical degree $10$ ample polarization $\sho_X(1)$, and let $\shu_X$ be the pullback of the tautological bundle on $\Gr(2,V_5)$ along the Gushel map $\gamma\colon X\rightarrow\Gr(2,V_5)$, then $\DCoh(X)$ has a geometric triangulated subcategory
$$\Ku(X):=\{E\in\DCoh(X)\,|\,\textsf{R}\Hom(\sho_X,E)=\textsf{R}\Hom(\shu^\vee_X,E)=0\}$$
which is called the \emph{Kuznetsov component} of $X$. 

The numerical Grothendieck group $\Num(\Ku(X))$ has rank 2. Denote the first Chern class of the polarization $\sho_X(1)$ by $h$, then one has a basis $\kappa_1,\kappa_2$ of $\Num(\Ku(X))$ with
$$\ch(\kappa_1)=-1+\frac{1}{5}h^2\quad\textup{and}\quad\ch(\kappa_2)=2-h+\frac{1}{12}h^3$$
such that the Euler form is given in this basis by $-I_2$. One consults \cite{KP17,PPZ26} for details.

Suppose now that $X$ is special, then the opposite of $X$ is a strongly smooth Gushel--Mukai surface $S$ and one has a covering involution $f\colon X\rightarrow X$ for the double cover $X\rightarrow M_S$ branched over the divisor $S\subset M_S$. The direct image of the covering involution restricts to an involution on the Kuznetsov component $\Ku(X)$ and induces a $\mathfrak{S}_2$-action on $\Ku(X)$. 

According to \cite{KP17,PPZ26}, the equivariant category $\Ku(X)_{\mathfrak{S}_2}$ is equivalent to $\DCoh(S)$ and the residual $\mathfrak{S}^\vee_2$-group action on $\DCoh(S)$ is induced by an involution on $\DCoh(S)$ say $\Pi$. The involution induces a Hodge isometry on $\Mukai(S,\bz)$ whose invariant lattice is spanned by the orthogonal Mukai vectors $v_1=(1,0,-1)$ and $v_2=(-2,H,-2)$.

The forgetful functor $\funct{Forg}\colon\DCoh(S)\rightarrow\Ku(X)$ induces a group homomorphism $\forg$ from the lattice $\Mukai(S,\bz)$ to the lattice $\Num(\Ku(X))$ and the inflation functor $\funct{Inf}\colon\Ku(X)\rightarrow\DCoh(S)$ induces a group homomorphism $\sfinf\colon \Num(\Ku(X))\rightarrow \Mukai(S,\bz)$.

\begin{proposition}[{{\cite[Proposition 6.6]{Liu26}}}]\label{En_image-forgetful-and-inflation}
In the setups above, one has
\begin{align*}
\forg(r,\Delta,s)&=(r-s)\kappa_1+(2r+2s+\Delta\cdot H)\kappa_2\\
\sfinf(a\kappa_1+b\kappa_2)&=av_1+bv_2
\end{align*}
for any vector $(r,\Delta,s)\in\Mukai(S,\bz)$ and integers $a,b$.
\end{proposition}

\begin{remark}
Here we choose $v_1$ and $v_2$ differently from \cite{Liu26}, because we want the semistable objects considered later to be in the heart. It is in fact not necessary.
\end{remark}

Moreover, the stability condition $\sigma_S$ in Example \ref{K3-Stab_the-stability-condition-on-GM-surface} is invariant under the involution $\Pi$ according to \cite[Proposition 6.2]{Liu26}.

\section{Double dual EPW sextics as Bridgeland moduli spaces}
In this section, we are going to show that 
$$\modulis_{\sigma_S}(v_1)\stackrel{\sim}{\rightarrow}\tilde{Y}^{\geq 1}_{A(S)^{\perp}}$$
for any strongly smooth Gushel--Mukai surface $S$. It is exactly Theorem~\ref{00_main-dual-sextic}.

\subsection{Birational model of the double dual EPW sextics}\label{subsection_birational_model}
To achieve the isomorphism, one needs the description of the associated double dual EPW sextic in \cite[Section 7.5]{DK24}.

One notices first that each line $L\subset S$ gives a Lagrangian submanifold $\bp^2\cong(\bp^1)^{[2]}$ in the hyperkähler manifold $S^{[2]}$. Then, one has a Mukai flop
\begin{displaymath}
	\xymatrix@R=2mm{
		&\hat{X}_S\ar[dr]\ar[dl]&\\
		S^{[2]}\ar[dr]\ar@{-->}[rr]^g&&X_S\ar[dl]\\
		&\bar{X}_S&
	}
\end{displaymath}
of these Lagrangian planes in $S^{[2]}$. The variety $\bar{X}_S$ is a subvariety of $\Gr(2,W)$ and the contraction $S^{[2]}\rightarrow \bar{X}_S$ is given by sending a point in $S^{[2]}$ to the line it spans in $\Gr(2,W)$.

Then one considers the Fano variety of lines $F(M_S)\cong\bp^2$ in the Grassmannian hull $M_S$ of the strongly smooth Gushel--Mukai surface $S$. It canonically contains the Fano variety of lines $F(S)$ on $S$, and embeds into the contracted space $\bar{X}_S$. So the blow-up $Z_S$ of $F(M_S)$ along the set $F(S)$ of points embeds into $S^{[2]}$. The strict transformation of the blow-up $Z_S$ along the Mukai flop $g$ is isomorphic to $F(M_S)$ and will be again denoted by $F(M_S)$. 

Also, one notices that a conic (including the degenerate ones) and a line in $S$ cannot meet exactly at a length $2$ subscheme due to Theorem~\ref{GM_nice-surface-condition}. So the dual plane $\Sigma_C^\vee$ parameterizing lines in the plane $\Sigma_C$ spanned by $C$ can be realized as a subvariety in $X_S$.

Indeed, suppose that $C$ is smooth, then every line in $\Sigma_C$ intersects $C$ in a length $2$ subscheme which is not contained in the exceptional locus of $g$ as the conic $C$ cannot meet a line in $S$ at a subscheme of length $2$. So $\Sigma_C^\vee\cong\Sigma_C^{[2]}\subset S^{[2]}$ is not contained in the exceptional locus of $g$ and hence embeds into $X_S$. The singular case requires a separate argument. Suppose for example that $C$ is a union of two lines, then $\Sigma_C^{[2]}$ is the union of the two Lagrangian planes spanned by these lines and the blow-up of $\Sigma_C^\vee$ along the two points represented by the lines. So one can see that $\Sigma_C^\vee$ is the strict transformation of $\Sigma_C^{[2]}$ via the flop $g$ and hence embeds into $X_S$. \footnote{These details are not included in \cite{DK24} and we thank Sasha Kuznetsov for the further explanation.}

The singular locus of $\tilde{Y}_{A(S)^{\perp}}^{\geq 1}$ is equal to $f_{A(S)^{\perp}}^{-1}(Y^{\geq 3}_{A(S)^{\perp}})$ as a set and for any $p$ in $Y^{\geq 3}_{A(S)^{\perp}}$ one uses $\tilde{p}$ to denote its preimage in $f_{A(S)^{\perp}}^{-1}(Y^{\geq 3}_{A(S)^{\perp}})$. One always has a distinguished point
$$p_S\in Y^{\geq 3}_{A(S)^{\perp}}=\{[v]\in\bp(V_6^\vee)\,|\,\dim(A(S)^{\perp}\cap\left(v\land\wedgep^2V_6^\vee\right))\geq 3\}$$
determined by the subspace $V_5\subset V_6$ according to \cite[Theorem 3.10]{DK18}, and $Y^{\geq 3}_{A(S)^{\perp}}-\{p_S\}$ is equal to the set of conics in $S$. The point $p_C\in Y^{\geq 3}_{A(S)^{\perp}}$ corresponding to a conic $C\subset S$ is given by the unique five-dimensional subspace $U_5\subset V_6$ determined by the linear span of $C$ in $\bp(W)$.

\begin{theorem}[{{\cite[Theorem 7.7 and Remark 7.9]{DK24}}}]\label{EPW-sextics-dual_geometry}
There is a diagram 
$$S^{[2]}\stackrel{g}{\dashrightarrow} X_S\stackrel{f^+}{\longrightarrow}\tilde{Y}_{A(S)^{\perp}}^{\geq 1}$$
where the birational map $g$ is the Mukai flop of all Lagrangian planes in $S^{[2]}$ corresponding to lines in $S$ and the morphism $f^+$ is a symplectic resolution of singularities. In particular, the map $f^+$ contracts $F(M_S)\subset X_S$ to the singular point $\tilde{p}_S$ and, for any other point corresponding to a conic $C\subset S$, the dual projective plane $\Sigma^\vee_C$ to that point.
\end{theorem}

\begin{remark}
The Mukai flop $g$ coincides with the flop in \cite[Equation (4.2.4)]{O'Grady13}. Once $S$ contains no lines or conics, the resolution $f^+$ coincides with the one in \cite[Equation (4.2.2)]{O'Grady13}.
\end{remark}

In the following subsections, we will recover the diagram by studying the wall-crossing behavior of the moduli space $\modulis_{\sigma_t}(v_1)$ for the stability condition $\sigma_t:=\sigma_{\sqrt{t}/5H,-2/5H}$ with $t\geq 1$. 

\subsection{Exclusion of the non-flopping walls}
In this subsection, we will show that either $\sigma_t$ is not on a wall with respect to $v_1$ or is on a flopping wall with respect to $v_1$.

The central charge of the stability condition $\sigma_t$ is 
$$Z_{\sigma_t}(r,\Delta,s)=\left((\frac{t-4}{5})r-\frac{2}{5}H\cdot\Delta-s\right)+\left(\frac{4}{5}r+\frac{1}{5}\Delta\cdot H\right)\sqrt{-t}$$
and the associated lattice $\bbmH_W$ of a potential wall $W$ for $v_1=(1,0,-1)$ is characterized by
$$\Im\frac{Z_{\sigma_t}(r,\Delta,s)}{Z_{\sigma_t}(1,0,-1)}=\Im\left(\frac{((t-4)r-2\Delta\cdot H-5s)+(\Delta\cdot H+4r)\sqrt{-t}}{(t+1)+4\sqrt{-t}}\right)=0$$
or equivalently the condition $20r+(t+9)\Delta\cdot H+20s=0$.

\begin{proposition}\label{EPW-sextics-dual_wall-crossing-BN}
There is no class $w\in\bbmH_W$ such that $\langle w,w\rangle=-2$ and $\langle v,w\rangle=0$.

\begin{proof}
Suppose otherwise that there is such a vector $w=(r,\Delta,s)$, then one has $(t+9)\Delta\cdot H=-40r$ and $\Delta^2=2r^2-2$. It follows
$$20r^2-20=10\Delta^2\leq (\Delta\cdot H)^2=\left(\frac{40r}{t+9}\right)^2$$
or equivalently $((t+9)^2-80)r^2\leq(t+9)^2$. So $r^2\leq5$ and all the possible cases are
\begin{itemize}
	\item $r=0$, $t\geq1$ and $\Delta^2=-2$ and $\Delta\cdot H=0$;
	\item $r=\pm1$, $t=1$, $\Delta^2=0$ and $\Delta\cdot H=\pm4$;
	\item $r=\pm1$, $t=13/3$, $\Delta^2=0$ and $\Delta\cdot H=\pm3$;
	\item $r=\pm1$, $t=11$, $\Delta^2=0$ and $\Delta\cdot H=\pm2$;
	\item $r=\pm1$, $t=31$, $\Delta^2=0$ and $\Delta\cdot H=\pm1$;
	\item $r=\pm2$, $t=1$, $\Delta^2=6$ and $\Delta\cdot H=\pm8$.
\end{itemize}
The first case is absurd as $|\Delta|$ should contain an effective class and $H$ is very ample, the third case will contradict Example \ref{GM_surface-condition}, and the other cases contradict Proposition~\ref{GM_nice-surface-condition} (in the last case one has $(H\pm\Delta)^2=0$ and $(H\pm\Delta)\cdot H=2$). 
\end{proof}
\end{proposition}

\begin{proposition}\label{EPW-sextics-dual_wall-crossing-HC}
There is no vector $w\in\bbmH_W$ such that $\langle w,w\rangle=0$ and $\langle v,w\rangle=1$.

\begin{proof}
Suppose otherwise that there exists such a vector $w=(r,\Delta,s)$, then one has $(t+9)\Delta\cdot H=-40r+20$ and $\Delta^2=2r^2-2r$. It follows
$$20r^2-20r=10\Delta^2\leq (\Delta\cdot H)^2=\left(\frac{40r-20}{t+9}\right)^2$$
or equivalently $((t+9)^2-80)r^2-((t+9)^2-80)r\leq20$. So $r\in\{0,1\}$ and all possible cases are
\begin{itemize}
	\item $r=0,t=1$ and $\Delta^2=0$ and $\Delta\cdot H=2$;
	\item $r=0,t=11$ and $\Delta^2=0$ and $\Delta\cdot H=1$; 
	\item $r=1,t=1$ and $\Delta^2=0$ and $\Delta\cdot H=-2$; 
	\item $r=1,t=11$ and $\Delta^2=0$ and $\Delta\cdot H=-1$; 
\end{itemize}
which always conflict with Proposition~\ref{GM_nice-surface-condition}.
\end{proof}
\end{proposition}

\begin{proposition}
There is no vector $w\in\bbmH_W$ such that $\langle w,w\rangle=0$ and $\langle v,w\rangle=2$.

\begin{proof}
Suppose otherwise that there exists such a vector $w=(r,\Delta,s)$, then one
has $(t+9)\Delta\cdot H=-40r+40$ and $\Delta^2=2r^2-4r$. It follows
$$20r^2-40r=10\Delta^2\leq (\Delta\cdot H)^2=\left(\frac{40r-40}{t+9}\right)^2$$
or equivalently $((t+9)^2-80)r^2-2((t+9)^2-80)r\leq80$. So $r\in\{-1,0,1,2,3\}$ and one can check as in Proposition~\ref{EPW-sextics-dual_wall-crossing-BN} and Proposition~\ref{EPW-sextics-dual_wall-crossing-HC} that all possible cases are invalid.
\end{proof}
\end{proposition}

\begin{proposition}
The potential wall cannot be a totally semistable wall.

\begin{proof}
Thanks to Proposition~\ref{EPW-sextics-dual_wall-crossing-HC}, it suffices to show that there is no $w=(r,\Delta,s)$ in $\bbmH_W$ with $\langle w,w\rangle=-2$ satisfying $r-s<0$ and $\Re(Z_{\sigma_t}(r,\Delta,s)/Z_{\sigma_t}(1,0,-1))>0$. Otherwise one has
\begin{align*}
r&=-\frac{t+9}{40}\Delta\cdot H-\frac{1}{2}\sqrt{\frac{(t+9)^2}{400}(\Delta\cdot H)^2-4-2\Delta^2}\\
s&=-\frac{t+9}{40}\Delta\cdot H+\frac{1}{2}\sqrt{\frac{(t+9)^2}{400}(\Delta\cdot H)^2-4-2\Delta^2}
\end{align*}
since $2rs=\Delta^2+2$ and $20(r+s)=-(t+9)\Delta\cdot H$.

The condition $\Re(Z_{\sigma_t}(r,\Delta,s)/Z_{\sigma_t}(1,0,-1))>0$ is equivalent to
$$(t^2+13t-4)r+(2t-2) \Delta\cdot H -5(t+1)s>0$$
and, after plugging $s$ and $r$, it becomes
$$(t^2+18t+1)\left(r+\frac{1}{4}\Delta\cdot H\right)>0$$
which implies $\Delta\cdot H+4r>0$ as $t\geq 1$. Together with $40r<-(t+9)\Delta\cdot H$ and $t\geq 1$, one can see that $\Delta\cdot H<0$ and $r>0$. Also, since $r<s$ one has $s\geq r+1$ so that
$$16r^2>(\Delta\cdot H)^2\geq 10\Delta^2=20rs-20\geq 20r(r+1)-20$$
which is not possible for integer $r>0$.
\end{proof}
\end{proposition}

In summary, every actual wall crossed by the path $\sigma_t$ is a flopping wall.

\subsection{The wall-crossing of spherical factors}
In this subsection, some stable spherical objects will be identified as the Jordan--Hölder factors appearing in the next subsection.

\begin{proposition}
Consider a line $L\subset S$, then $\modulis_{\sigma_t}(0,L,-1)=\modulis_H(0,L,-1)$ for any $t\geq 1$.
	
\begin{proof}
It suffices to show that $\modulis_{\sigma_t}(0,L,-1)$ contains only stable objects in the heart of $\sigma_t$ for any $t\geq 1$. Otherwise, one has an exact sequence $0\rightarrow A\rightarrow \sho_L(-2)\rightarrow B\rightarrow0$ of $\sigma_t$-semistable objects with the same phase. Set $v(A)=(r,\Delta,s)$, then $r\geq 1$ and one has
$$\Im(Z_t(A))=(\frac{1}{5}\Delta\cdot H+\frac{4}{5}r)\sqrt{t}>0\quad\textup{and}\quad \Im(Z_t(B))=(\frac{1}{5}-\frac{1}{5}\Delta\cdot H-\frac{4}{5}r)\sqrt{t}>0$$
implying that $\Delta\cdot H+4r$ cannot be an integer, which is absurd.
\end{proof}
\end{proposition}

\begin{proposition}
Consider a line $L\subset S$, then $\modulis_{\sigma_t}(1,-L,0)=\modulis_H(1,-L,0)$ for any $t\geq 1$.
	
\begin{proof}
It suffices to show that $\modulis_{\sigma_t}(1,-L,0)$ contains only stable objects in the heart of $\sigma_t$ for any $t\geq 1$. Otherwise, one has an exact sequence $0\rightarrow A\rightarrow \sho_S(-L)\rightarrow B\rightarrow0$ of $\sigma_t$-semistable objects with the same phase such that $A$ is stable. Set $v(A)=(r,\Delta,s)$, then 
$$\Im(Z_t(A))=(\frac{1}{5}\Delta\cdot H+\frac{4}{5}r)\sqrt{t}>0\quad\textup{and}\quad \Im(Z_t(B))=(\frac{3}{5}-\frac{1}{5}\Delta\cdot H-\frac{4}{5}r)\sqrt{t}>0$$
so that $k:=\Delta\cdot H+4r$ is either $1$ or $2$. It follows $3\Re(Z_t(A))=k\Re(Z_t(\sho_S(-L)))$ and then an equation $15s=(3t+12)r-4k-kt$.

Since $(1,-L,0)$ is a spherical vector, the stable HN factor $A$ of $\sho_S(-L)$ is also spherical due to Proposition~\ref{K3-Stab_wall-criteria-spherical} (or directly \cite[Lemma 12.2]{Bri08}). So the Hodge index theorem gives
$$(4r-k)^2\geq10\Delta^2=10(2rs-2)=\frac{4}{3}((3t+12)r-4k-kt)r-20$$
and then the inequality $60+3k^2 \geq t(12r-4k)r+8kr\geq 12r^2+4kr$ as $3r\geq 3>k$. It follows that $r=1$ or $r=2$. The latter case is impossible as one can always see $s\geq 2$ from $t\geq 1$ and the equality for $s$. So the Hodge index theorem becomes $(8-k)^2\geq 10\Delta^2\geq 60$, impossible.

Hence one must have $r=1$, then $\shh^{-1}(B)=0$ from the destabilizing sequence. In this case, the object $B$ is also a spherical stable factor of $\sho_S(-L)$ by Proposition~\ref{K3-Stab_wall-criteria-spherical}. It follows that $s=1-\Delta\cdot L$. One notices that $\Delta\cdot H\leq -2$ as $k=1,2$. Then one can see from the Hodge index theorem that $s=1$, hence $\Delta^2=0$. So one has $\Delta\cdot H=-2$ or $\Delta\cdot H=-3$ which will contradict either Example \ref{GM_surface-condition} or Proposition~\ref{GM_nice-surface-condition}.
\end{proof}
\end{proposition}

\begin{proposition}
Consider a conic $C\subset S$, then $\modulis_{\sigma_t}(0,C,-1)=\modulis_H(0,C,-1)$ for $t\geq 1$.
	
\begin{proof}
It suffices to show that $\modulis_{\sigma_t}(0,C,-1)$ contains only stable objects in the heart of $\sigma_t$ for any $t\geq 1$. Otherwise, one has an exact sequence $0\rightarrow A\rightarrow \sho_C(-2)\rightarrow B\rightarrow0$ of $\sigma_t$-semistable objects with the same phase such that $A$ is $\sigma_t$-stable. Set $v(A)=(r,\Delta,s)$, then one has
$$\Im(Z_t(A))=(\frac{1}{5}\Delta\cdot H+\frac{4}{5}r)\sqrt{t}>0\quad\textup{and}\quad \Im(Z_t(B))=(\frac{2}{5}-\frac{1}{5}\Delta\cdot H-\frac{4}{5}r)\sqrt{t}>0$$
so that $\Delta\cdot H+4r=1$. Then $2\Re(Z_t(A))=\Re(Z_t(\sho_C(-2)))$ and one gets $10s=(2t+8)r-5$. 

Since $(0,C,-1)$ is a spherical vector, the stable factor $A$ is also spherical. So the Hodge index theorem gives $(4r-1)^2\geq10\Delta^2=10(2rs-2)=(4t+16)r^2-10r-20$ from which one obtains $r\in\{0,1,2\}$ as $r\geq 0$ and $t\geq 1$. Suppose that $r=0$, then one has $2s=-1$, which is absurd. Suppose that $r=1$, then $\Delta\cdot H=-3$. However, one can see from $10s=2t+3$ and the Hodge index theorem that $\Delta^2=0$, which is impossible as $(S,H)$ is Brill--Noether general of degree $10$. Suppose that $r=2$, then $s\geq 2$ and $\Delta^2\geq 6$, violating the Hodge index theorem. 
\end{proof}
\end{proposition}

\begin{proposition}
Consider a conic $C\subset S$, then $\modulis_{\sigma_t}(1,-C,0)=\modulis_H(1,-C,0)$ for $t\geq 1$.
	
\begin{proof}
It suffices to show that $\modulis_{\sigma_t}(1,-C,0)$ contains only stable objects in the heart of $\sigma_t$ for any $t\geq 1$. Otherwise, one has an exact sequence $0\rightarrow A\rightarrow \sho_S(-C)\rightarrow B\rightarrow0$ of $\sigma_t$-semistable objects with the same phase such that $A$ is $\sigma_t$-stable. Set $v(A)=(r,\Delta,s)$, then one has
$$\Im(Z_t(A))=(\frac{1}{5}\Delta\cdot H+\frac{4}{5}r)\sqrt{t}>0\quad\textup{and}\quad \Im(Z_t(B))=(\frac{2}{5}-\frac{1}{5}\Delta\cdot H-\frac{4}{5}r)\sqrt{t}>0$$
so that $\Delta\cdot H+4r=1$. Then $\Re(Z_t(A))=\Re(Z_t(B))$ and one gets $10s=(2t+8)r-t-4$. 

Since $(1,-C,0)$ is a spherical vector, the stable factor $A$ is also spherical. So the Hodge index theorem gives $(4r-1)^2\geq10\Delta^2=10(2rs-2)=(4t+16)r^2-(2t+8)r-20$ from which one obtains $r=1$ or $r=2$ as $r\geq 1$ and $t\geq 1$. Suppose that $r=1$, then $\Delta\cdot H=-3$. On the other hand, the Hodge index theorem and $10s=t+4>0$ ensure that $\Delta^2=0$. It is impossible since $(S,H)$ is a Brill--Noether general K3 surface. Suppose that $r=2$, then one can see $s\geq 1$ and $\Delta^2\geq 6$ but it will violate the Hodge index theorem. 
\end{proof}
\end{proposition}

Now one considers the unique Gieseker stable bundle $\shu_S$ on $S$ with Mukai vector $(2,-H,3)$ which is known to be the restriction of the tautological bundle $\shu_G$ on $\Gr(2,V_5)$.

\begin{proposition}\label{EPW-sextics-dual_flop3-factor1}
One has $\modulis_{\sigma_t}(-2,H,-3)=\modulis_H(2,-H,3)[1]$ for any $t\geq 1$.
	
\begin{proof}
It suffices to show that $\modulis_{\sigma_t}(-2,H,-3)$ contains only stable objects in the heart of $\sigma_t$ for any $t\geq 1$. Otherwise, one has an exact sequence $0\rightarrow A\rightarrow \shu_S[1]\rightarrow B\rightarrow0$ of $\sigma_t$-semistable objects with the same phase such that $B$ is $\sigma_t$-stable. Set $v(B)=(r,\Delta,s)$, then one sees from the destabilizing exact sequence that $r\leq-1$ as $\shh^{-1}(B)$ is torsion-free. Moreover, one can see as before that $\Delta\cdot H+4r=1$. It follows $\Re(Z_t(B))=\Re(Z_t(A))$ and $10s=(2t+8)r+2t-7$. 

In this case, both $A$ and $B$ are spherical stable objects as $\shu_S[1]$ is spherical. So one can get $2s=r-2$. Then $r=-1$ by the previous equation, which is impossible. 
\end{proof}
\end{proposition}

Then one considers the unique Gieseker stable bundle $\shv_S$ on $S$ with Mukai vector $(3,-H,2)$ which is actually the restriction of the tautological bundle $\shv_G$ on $\Gr(3,V_5)\cong\Gr(2,V_5)$.

\begin{proposition}\label{EPW-sextics-dual_flop3-factor2}
One has $\modulis_{\sigma_t}(3,-H,2)=\modulis_H(3,-H,2)$ for any $t\geq 1$.
	
\begin{proof}
It suffices to show that $\modulis_{\sigma_t}(3,-H,2)$ contains only stable objects in the heart of $\sigma_t$ for any $t\geq 1$. Otherwise, one has an exact sequence $0\rightarrow A\rightarrow \shv_S\rightarrow B\rightarrow0$ of $\sigma_t$-semistable objects with the same phase such that $A$ is $\sigma_t$-stable. Set $v(A)=(r,\Delta,s)$, then one can see $\Delta\cdot H+4r=1$. It follows $\Re(Z_t(A))=\Re(Z_t(B))$ and then $10s=(2t+8)r-2-3t$.

In this case, both $A$ and $B$ are spherical stable objects as $\shv_S$ is spherical. So one can get $3s=2r$. Then $s=1$ by the previous equation, which is impossible. 
\end{proof}
\end{proposition}

\subsection{The flopping walls and identification}
In this subsection, we will recover the birational map by identifying the possible flopping walls for the moduli space $\modulis_{\sigma_t}(v_1)$.

\begin{proposition}\label{EPW-sextics-dual_the-first-flop}
Suppose that $t>1$, then the stability condition $\sigma_t$ is on a flopping wall for $v_1=(1,0,-1)$ if and only if $t=11$ and $S$ contains lines. Moreover, the induced flop map
$$S^{[2]}\dashrightarrow Y_S$$
coincides with the Mukai flop $g$ in Theorem~\ref{EPW-sextics-dual_geometry}.

\begin{proof}
According to Theorem~\ref{K3-Stab_wall-criteria-flop}, a flopping wall for $v_1$ can only be induced by a pair of vectors $w$ and $w'$ in $\bbmH_W$ such that $\langle w,w\rangle=\langle w',w'\rangle=-2,\langle v_1,w\rangle=\langle v_1,w'\rangle=1$ and $w'+w=v_1$. One writes $w=(r,\Delta,s)$, then $(t+9)\Delta\cdot H=-40r+20$ and $\Delta^2=2r^2-2r-2$. It follows
$$20r^2-20r-20=10\Delta^2\leq (\Delta\cdot H)^2=\left(\frac{40r-20}{t+9}\right)^2$$
or equivalently $r^2-r-1\leq100/((t+9)^2-80)$. It follows $t=3,11$ as $t>0$ and $r\in\bz$.

Suppose that $t=3$, then $r\in\{-1,2\}$. It means that $\Delta^2=2$ and $\Delta\cdot H=\pm 5$. However, the class $2\Delta\pm H$ will have square $-2$ and degree zero. It is impossible as $H$ is ample.

Hence $t=11$ and $r\in\{1,0\}$. Then all possible vectors are $(1,\Delta,0)$ with $\Delta^2=-2,\Delta\cdot H=-1$ and $(0,\Delta,-1)$ with $\Delta^2=-2,\Delta\cdot H=1$.

The divisor class $\Delta$ in $(0,\Delta,-1)$ is in fact the class of a line in $S$ (see, for example, the discussions in \cite[Section 2.1.4]{HuK3}) and the class $\Delta$ in $(1,\Delta,0)$ is $-L$ for some line $L$ in $S$.

According to Theorem~\ref{K3-Stab_K3-large-volume-limit}, one has $\modulis_{\sigma_t}(1,0,-1)=\modulis_H(1,0,-1)\cong S^{[2]}$ for $t>11$. Then the flopping walls for $\modulis_{\sigma_{11}}(1,0,-1)$ are induced by the short exact sequences
$$0\rightarrow\sho_S(-L)\rightarrow\shi_Z\rightarrow \sho_L(-Z)\rightarrow0$$
where $Z$ is any length $2$ subscheme on any line $L\subset S$. Then one has a flop
$$g\colon S^{[2]}\dashrightarrow Y_S=\modulis_{\sigma_t}(1,0,-1)$$
for any $11>t>1$, which coincides with the Mukai flop in Theorem~\ref{EPW-sextics-dual_geometry} by construction. 
\end{proof}
\end{proposition}

It remains to show that $Y_S\cong X_S$ contracts to the moduli space $\modulis_{\sigma_S}(1,0,-1)$ through the same symplectic resolution $f^+\colon X_S\rightarrow\tilde{Y}_{A(S)^{\perp}}^{\geq 1}$ as in Theorem~\ref{EPW-sextics-dual_geometry}.

\begin{proposition}\label{EPW-sextics-dual_the-second-flop}
The moduli space $\modulis_{\sigma_S}(v_1)$ is isomorphic to $\tilde{Y}^{\geq 1}_{A(S)^{\perp}}$.

\begin{proof}
According to Theorem~\ref{K3-Stab_wall-criteria-flop}, a flopping wall for $\modulis_{\sigma_S}(v_1)$ is induced by a pair of vectors $w$ and $w'$ in $\bbmH_W$ such that $\langle w,w\rangle=\langle w',w'\rangle=-2,\langle v_1,w\rangle=\langle v_1,w'\rangle=1$ and $w'+w=v_1$. One writes $w=(r,\Delta,s)$, then $\Delta\cdot H=-4r+2$ and also $\Delta^2=2r^2-2r-2$. It follows
$$20r^2-20r-20=10\Delta^2\leq (\Delta\cdot H)^2=(4r-2)^2$$
or equivalently $r^2-r-1\leq 5$. So one has $r\in\{-2,-1,0,1,2,3\}$. 

At first, one notices that $r$ cannot be $-1$ or $2$ as in this case one will get a class $\Delta$ with $\Delta\cdot H=6$ and $\Delta^2=2$. It is impossible as otherwise $H-\Delta$ would violate Proposition~\ref{GM_nice-surface-condition}. 

Also, one can see that the class $\Delta$ in $w=(1,-\Delta,0)$ and $w'=(0,\Delta,-1)$ with $\Delta\cdot H=2$ and $\Delta^2=-2$ is in fact the class of a conic $C$ in $S$ (see for example \cite[Section 2.1.4]{HuK3}). The corresponding flopping wall is induced by the short exact sequences
$$0\rightarrow\sho_S(-C)\rightarrow\shi_Z\rightarrow \sho_C(-Z)\rightarrow0$$
for the conic $C\subset S$ and any length two subscheme $Z\subset C\subset S$ which is not in the exceptional locus of the Mukai flop $g\colon S^{[2]}\dashrightarrow X_S$. So the space $\Sigma_C^\vee$ is contracted to a point.

Then one considers the pair of vectors $w=(-2,\Delta,-3)$ and $w'=(3,-\Delta,2)$ with $\Delta\cdot H=10$ and $\Delta^2=10$. The Hodge index theorem ensures that $\Delta=H$. So according to Proposition~\ref{EPW-sextics-dual_flop3-factor1} and Proposition~\ref{EPW-sextics-dual_flop3-factor2}, the corresponding destabilizing exact sequence for $\shi_Z$ is
$$0\rightarrow\shv_S\rightarrow\shi_Z\rightarrow\shu_S[1]\rightarrow0$$
where $Z\subset S$ is a length $2$ subscheme spanning a line not in $S$. Indeed, such a short exact sequence follows from the fact that the ideal sheaf $\shi_{L/M_S}$ of every line $L\subset M_S$ can be written as the cokernel of a unique map $\shu_G|_{M_S}\rightarrow\shv_G|_{M_S}$, and each such map has the ideal sheaf of a line as the cokernel. One just takes the restriction onto $S\subset M_S$. The line spanned by $Z$ is not in $S$ because such a length two subscheme $Z$ has already been flopped in Proposition~\ref{EPW-sextics-dual_the-first-flop}. So one contracts the Fano variety $F(M_S)$ of lines, and therefore recovers $f^+\colon X_S\rightarrow \tilde{Y}_{A(S)^{\perp}}^{\geq 1}$.

Moreover, the computations above show that the contracted space parameterizes all the $S$-equivalence classes of $\sigma_S$-semistable objects with class $v_1$. This proves the claim.
\end{proof}
\end{proposition}

In conclusion, we finish the proof of Theorem~\ref{00_main-dual-sextic}.

\begin{remark}\label{EPW-sextics-dual_the-involution}
Moreover, the involution on $\modulis_{\sigma_S}(v_1)$ induced by the autoequivalence $\Pi$ coincides with the covering involution on the double dual EPW sextic as the Hodge isometries of the corresponding birational involutions on $X_S$ are the same (see \cite[Section 6.5]{Liu26} for the Picard rank one case). A point-wise description of the involution can be found in \cite[Remark 7.8]{DK24}, and the fixed locus $\modulis_{\sigma_S}(v_1)_{\mathfrak{S}^\vee_2}$ of the involution is exactly the ramification locus $Y^{\geq 2}_{A(S)^{\perp}}$.
\end{remark}

\section{Other EPW varieties as Bridgeland moduli spaces}
In this section, we are going to realize other EPW varieties as Bridgeland moduli spaces and in particular prove Proposition~\ref{00_main-sextic-weak} and Proposition~\ref{00_main-dual-surface}.

\subsection{Double EPW sextics}
It is shown in \cite[Section 6.6]{Liu26} that the moduli space $\modulis_{\sigma_S}(v_2)$ is isomorphic to the double EPW sextic once $S$ has Picard rank one. Here we extend the result to the case where $Y_{A(S)}^{\geq 3}=\varnothing$. According to \cite[Theorem 3.3]{Beri26}, it means that $S$ contains no lines or quintic elliptic curves, allowing us to achieve the following two observations.
 
\begin{proposition}
Suppose that $(S,H)$ is a strongly smooth Gushel--Mukai surface containing no quintic elliptic curves, then $\modulis_H(-v_2)$ is a hyperkähler manifold.
	
\begin{proof}
It suffices to show that there are no strictly semistable sheaves in $\modulis_H(-v_2)$. Suppose, by contradiction, that there exists such a strictly semistable sheaf $\shg$. Then one can take a Jordan--Hölder factor $\shb$ of $\shg$ with $v(\shb)=(1,\Delta,s)$. One can see
$$\Delta\cdot H=5\quad\textup{and}\quad s=1$$
by the equality of reduced Hilbert polynomials. So $\shb\cong\shi_Z(\Delta)$ for some zero-dimensional subscheme $Z\subset S$. It follows that $\Delta^2\geq0$ as the length of $Z$ is non-negative.
		
By the Hodge index theorem, one can see that $\Delta^2=0$ or $\Delta^2=2$. The case $\Delta^2=0$ will violate our assumption, while the case $\Delta^2=2$ will violate Example \ref{GM_surface-condition}.
\end{proof}
\end{proposition}

\begin{proposition}\label{EPW-sextics_generic-no-walls}
Suppose that $(S,H)$ is a strongly smooth Gushel--Mukai surface containing no lines or quintic elliptic curves, then one has $\modulis_{\sigma_S}(v_2)\cong \modulis_H(-v_2)$.

\begin{proof}
It suffices to show that $\sigma_t$ is $v_2$-generic for any $t\geq 1$. Otherwise, one chooses the largest $t\geq 1$ such that $\sigma_t$ is on a wall with respect to $v_2$. Then a strictly $\sigma_t$-semistable object in the heart with Mukai vector $v_2$ has the form $\shg[1]$ for some $\shg$ in $\modulis_H(-v_2)=\modulis_H(2,-H,2)$.

One finds a short exact sequence $0\rightarrow A\rightarrow \shg[1]\rightarrow B\rightarrow0$ of $\sigma_t$-semistable objects with the same phase such that $A$ is stable. One can first see $r\geq-1$ for $v(A)=(r,\Delta,s)$ from the long exact sequence of cohomology sheaves $0\rightarrow \shh^{-1}(A)\rightarrow \shg\rightarrow \shh^{-1}(B)\rightarrow\shh^0(A)\rightarrow0$.

Also, one can obtain
$$
\Delta\cdot H=1-4r\quad\textup{and}\quad s=\frac{t+4}{5}r+\frac{t-1}{5}
$$
through $\mu(A)=\mu(B)=\mu(\shg[1])<+\infty$. Then one can see
$$(4r-1)^2=(\Delta\cdot H)^2\geq10\Delta^2\geq 20rs-20=4(t+4)r^2+4(t-1)r-20$$
using the Hodge index theorem and Mukai lemma. It follows
$$\frac{25}{4}\geq (tr+1)(r+1)$$
and one can see $r\leq 1$ from this inequality and $t\geq 1$. 

Suppose that $r=1$, then $\Delta\cdot H=-3$ and $9\geq 10\Delta^2\geq 8t-8\geq 8$. It follows $\Delta^2=0$, but then the class $-\Delta$ will violate Example \ref{GM_surface-condition}. 

Suppose that $r=0$, then $\Delta\cdot H=1$ and $1\geq 10\Delta^2\geq -20$. It follows $\Delta^2=-2$ due to Proposition~\ref{GM_nice-surface-condition}, but then the class $\Delta$ is a line which is impossible. 

Suppose that $r=-1$, then $\Delta\cdot H=5$ and $25\geq 10\Delta^2\geq 0$. It follows $\Delta^2=2$, as $S$ also does not contain quintic elliptic curves. However, in this case the class $D:=H-2\Delta$ is a degree zero effective class, which is impossible as $H$ is ample.
\end{proof}
\end{proposition}

Using the Mukai isomorphism
$$\theta_{-v_2}\colon(-v_2)^{\perp}\stackrel{\sim}{\rightarrow} H^2(\modulis_H(-v_2),\bz)$$
in \cite{BM14,BM14-MMP}, one can define a square two class $\theta_S:=\theta_{-v_2}(v_1)$. According to \cite[Theorem 12.1]{BM14-MMP} and the wall-crossing criterion in Section~\ref{K3-Stab_BM-criterion}, one can verify that $\theta_S$ is ample by Proposition~\ref{EPW-sextics_generic-no-walls}.

\begin{proposition}
Suppose that $(S,H)$ is a strongly smooth Gushel--Mukai surface containing no lines or quintic elliptic curves, then one has a polarized isomorphism $$(\modulis_H(-v_2),\theta_S)\stackrel{\sim}{\rightarrow} (\tilde{Y}_{A(S)}^{\geq 1},h_{A(S)})$$
of hyperkähler manifolds where $h_{A(S)}$ is the ample class corresponding to $\sho_{\tilde{Y}_{A(S)}^{\geq 1}}(1)$.

\begin{proof}
The statement holds when $S$ has Picard rank one according to \cite[Section 5.4]{O'Grady05} and the discussion in \cite[Section 6.6]{Liu26}. Then we want to spread this polarized isomorphism. 
	
To save notation, let us set the hyperkähler manifolds
$$
M_S:=\modulis_H(-v_2)\quad\textup{and}\quad Y_S:=\tilde{Y}_{A(S)}^{\geq 1}
$$
for a strongly smooth Gushel--Mukai surface $S$ containing no lines or quintic elliptic curves.

Now pick any such Gushel--Mukai surface $(S,H)$, and consider a small irreducible analytic neighborhood $B$ of the point represented by $(S,H)$ in the moduli space of Gushel--Mukai surfaces such that any point in $B$ corresponds to a strongly smooth Gushel--Mukai surface containing no lines or quintic elliptic curves. After shrinking $B$, we may assume that there is a family
$$
(\schs,\schh)\rightarrow B
$$
whose fiber over $b\in B$ is $(S_b,H_b)$. The associated Lagrangian
subspaces $A(S_b)$ vary in families according to \cite{DK20b,DK20c}, and hence give a smooth projective
family
$$
\schy\rightarrow B
$$
together with a relative EPW polarization $h_B$ on $\schy$ such that $(\schy_b,h_b)=(\tilde{Y}^{\geq 1}_{A(S_b)},h_{A(S_b)})$.

Similarly, the relative moduli spaces (see \cite{HL}) parameterizing $H_b$-stable sheaves on $S_b$ with Mukai vector $(2,-H_b,2)$ form a smooth projective family
$$
\schm:=\schm_{\schh}(2,-\schh,2)\rightarrow B
$$
which carries relative polarization class $\theta_B$ such that $\theta_b:=\theta_{S_b}$ for any $b\in B$.

By construction, both families define holomorphic maps
$$
B\longrightarrow \mathfrak{N}
$$
to the same connected component $\mathfrak{N}$ of the moduli space of polarized hyperkähler manifolds of $K3^{[2]}$-type with square two
polarizations. These holomorphic maps agree on the analytically dense subset $B^{\circ}\subset B$ containing Picard rank one K3 surfaces. Since $\mathfrak{N}$ is quasi-projective (see, for example, \cite{DM19}), it is separated. Then the two holomorphic maps agree on all of $B$ as $B$ is irreducible. In particular, the polarized periods of $(M_S,\theta_S)$ and $(Y_S,h_{A(S)})$ agree.

In this case, the global Torelli theorem for hyperkähler manifolds \cite[Theorem 1.3]{Mar11} ensures the desired isomorphism.
\end{proof}
\end{proposition}

Using the two identifications and Remark \ref{EPW-sextics-dual_the-involution}, one finishes Proposition~\ref{00_main-sextic-weak}.

\begin{proposition}\label{EPW-sextics_generic}
Suppose that $(S,H)$ is a strongly smooth Gushel--Mukai surface containing no lines or quintic elliptic curves, then one has $\modulis_{\sigma_S}(v_2)\cong \tilde{Y}_{A(S)}^{\geq 1}$ such that the involution on $\modulis_{\sigma_S}(v_2)$ induced by $\Pi$ coincides with the covering involution on $\tilde{Y}_{A(S)}^{\geq 1}$.
\end{proposition}

\subsection{Double dual EPW surfaces}
The statement in this subsection is essentially a corollary of \cite[Theorem 7.3]{DK24}, which describes the double dual EPW surface associated with a Gushel--Mukai threefold, and \cite[Theorem 7.12]{JLLZ24} which describes the corresponding moduli space.

Here we only elaborate on the case for a special Gushel--Mukai threefold $X$ over $S$ and the starting point is the Hilbert scheme $C(X)$ of conics on $X$. The scheme $C(X)$ is a normal projective surface which has two components $C_{\sigma}(X)$ and $C_0(X)$. The component 
$$C_{\sigma}(X)\cong F(M_S)\cong\bp^2$$ 
parameterizes the preimage of lines on $M_S$. The other component $C_0(X)$ is isomorphic to the blow-up of the double dual EPW surface $\tilde{Y}^{\geq2}_{A(X)^{\perp}}$ associated with $X$ at the distinguished point $\tilde{p}_S$ mentioned in Theorem~\ref{EPW-sextics-dual_geometry} according to \cite[Theorem 7.3]{DK24}. The two components intersect transversely along the exceptional curve on $C_0(X)$.

By \cite[Theorem 7.12]{JLLZ24}, the moduli space $\modulis_{\sigma_X}(\kappa_1)$ is a contraction from $C(X)$ such that the component $C_{\sigma}(X)$ is contracted to the distinguished singular point $\tilde{p}_S$. It follows that the moduli space $\modulis_{\sigma_X}(\kappa_1)$ is exactly the double dual EPW surface associated with $X$.

In conclusion, we can finish Corollary~\ref{00_main-dual-surface}.

\begin{proposition}\label{EPW-surface-dual_special}
Consider a strongly smooth Gushel--Mukai surface $S$ and the special Gushel--Mukai threefold $X$ over it, then $\modulis_{\sigma_X}(\kappa_1)$ is isomorphic to $\tilde{Y}^{\geq 2}_{A(S)^{\perp}}$ and the branched cover
$$\modulis_{\sigma_X}(\kappa_1)\rightarrow \modulis_{\sigma_S}(v_1)_{\mathfrak{S}^\vee_2}$$
induced by the forgetful functor coincides with $\tilde{Y}^{\geq 2}_{A(S)^{\perp}}\rightarrow Y^{\geq 2}_{A(S)^{\perp}}$.

\begin{proof}
The isomorphism $\modulis_{\sigma_X}(\kappa_1)\cong\tilde{Y}^{\geq 2}_{A(S)^{\perp}}$ follows from \cite{DK24,JLLZ24} as above. The identification of the branched double cover is similar to \cite[Proposition 5.20]{Liu26}. One can do this because the involution on $C(X)$ induced by the branched cover $X\rightarrow M_S$ coincides with the one induced by the natural covering involution on $\tilde{Y}^{\geq 2}_{A(S)^{\perp}}$ according to \cite[Section 7.3]{DK24}. In particular, the branched double cover for the double dual EPW surface $\tilde{Y}^{\geq 2}_{A(S)^{\perp}}$ can be described in terms of conics on an open subset and therefore be identified.
\end{proof}
\end{proposition}

The double cover can be described pointwise. The dual EPW surface $Y^{\geq 2}_{A(S)^{\perp}}$ is the union of a distinguished point $p_S$ and the surface $BTC_S$ of conics bitangent to $S$, where we adopt the convention that $BTC_S$ also parameterizes the conics in $S$. Any singular point $p_C$ in $ Y^{\geq 2}_{A(S)^{\perp}}$ other than $p_S$ corresponds to a unique conic $C$ contained in $S$.

The singular point in $\tilde{Y}^{\geq 2}_{A(S)^{\perp}}$ induced by the contraction from $C(X)$ is precisely $\tilde{p}_S$, and any other singular point corresponds to a conic contained in the ramification locus of the double cover $X\rightarrow M_S$. In this case, a smooth point in $\tilde{Y}^{\geq 2}_{A(S)^{\perp}}$ is represented by a conic bitangent to the ramification locus and corresponds to a conic in $M_S$ bitangent to $S$. Conversely, the preimage of a conic in $M_S$ bitangent to $S$ is the union of two distinct conics.

\section{The conjecture and implications}
\subsection{The duality conjecture}
Now let us explain our Conjecture \ref{00_conjecture}.

\begin{conjecture}\label{Conjecture}
Consider two Gushel--Mukai threefolds $X$ and $X'$ such that one has an isomorphism $\tilde{Y}_{A(X)}^{\geq 1}\cong\tilde{Y}_{A(X')^{\perp}}^{\geq 1}$ preserving the canonical polarizations, then there exists an equivalence $\cate{Ku}(X)\cong\cate{Ku}(X')$ such that $\kappa_1$ is sent to $\kappa'_2$ and $\kappa_2$ is sent to $\kappa'_1$ up to signs.
\end{conjecture}

It is in fact a refinement of the homological projective duality theorem in \cite{KP23}.

\begin{theorem}[{{\cite[Theorem 3.3.35]{KP23}}}]\label{Conjecture_BP-statement}
Consider two Gushel--Mukai threefolds $X$ and $X'$ such that there exists an isomorphism between $\tilde{Y}_{A(X)}^{\geq 1}$ and $\tilde{Y}_{A(X')^{\perp}}^{\geq 1}$ which preserves the canonical polarizations, then one has an equivalence $\cate{Ku}(X)\cong\cate{Ku}(X')$.
\end{theorem}

The equivalence is not canonically constructed in general. In fact, the theorem constructs a canonical equivalence using \cite{KP21} for special period duals, called the \emph{direct period duals}, and then uses Proposition~\ref{GM_threefold-period-dual} to conclude the general case. 

\begin{definition}
Let $X$ and $X'$ be two Gushel--Mukai threefolds such that $X'$ is a period dual of $X$. Then $X'$ is said to be a \defi{direct period dual} of $X$ if it corresponds to a point which is not in $\bp(V_5)\subset\bp(V_6)$ in the sense of Proposition~\ref{GM_threefold-period-dual}.
\end{definition}

The conjecture is true if $X$ is very general in the sense that it is not the period dual of itself. Hence the canonical equivalence $\Ku(X)\cong\Ku(X')$ constructed in Theorem~\ref{Conjecture_BP-statement} satisfies Conjecture \ref{Conjecture} for a very general direct period dual $X$ and $X'$.

\begin{proposition}\label{Conjecture-is-true-generically}
Conjecture \ref{Conjecture} holds if at least one of the two Gushel--Mukai threefolds $X$ and $X'$ is not the period dual of itself.
\end{proposition}

It is shown in \cite[Appendix C]{DK18} and \cite{DK20b} that the construction of Gushel--Mukai varieties can be generalized to the relative setting. In particular, one has the classical projective duality for relative direct period duals according to \cite[Lemma C.6]{DK18}.

Suppose that one has a relative version of Theorem~\ref{Conjecture_BP-statement} for a certain small base with positive dimension, then one can spread Proposition~\ref{Conjecture-is-true-generically} to any direct period duals. In this case, it is not difficult to verify Conjecture \ref{Conjecture}. This will appear in future work.

\subsection{Possible self-dual surfaces}
If there are no self-dual Gushel--Mukai threefolds, then it follows from Theorem~\ref{GM_period-double-EPW} that the equivalence $\Ku(X)\simeq\Ku(X')$ would send $\kappa_1$ to $\kappa_2'$ and $\kappa_2$ to $\kappa_1'$, up to sign. Unfortunately, this is not true by the example constructed in \cite{DM22}. We want to discuss the possible self-dual examples in our setting. 

Let $S$ be a strongly smooth Gushel--Mukai surface, and let $X$ be the special Gushel--Mukai threefold over it. Suppose that one has an autoequivalence on $\Ku(X)$ which exchanges $\kappa_1$ and $\kappa_2$ up to signs. Then, according to \cite[Lemma 4.9]{BP23}, the autoequivalence descends to an autoequivalence $\Psi$ on $\cate{D}^b(S)$ which exchanges $v_1$ and $v_2$ up to signs.

\begin{proposition}
Up to a shift on $\Psi$, the Hodge isometry $\psi$ on $\Mukai(S,\bz)$ induced by $\Psi$ satisfies $\psi(v_1)=v_2$ and $\psi(v_2)=-v_1$. Moreover, one of the following two cases holds:
\begin{itemize}
	\item[(a)] $\psi(0,0,1)=(0,E,-2)$ for a class $E$ satisfying $E^2=0$ and $E\cdot H=5$;
	\item[(b)] $\psi(0,0,1)=(1,L,-1)$ for a class $L$ satisfying $L^2=-2$ and $L\cdot H=1$.
\end{itemize}

\begin{proof}
One notices that $\psi$ cannot exchange $v_1$ and $v_2$ directly, otherwise it would be orientation-reversing which violates \cite{HMS09}. Then one computes as in \cite[Proof 2 in Page 142]{BP23}.
\end{proof}
\end{proposition}

Set $\sigma'_S:=\Phi.\sigma_S$, then one should have $\modulis_{\sigma_S}(0,0,1)\cong\modulis_{\sigma'_S}(\psi(0,0,1))$. Suppose further that the stability condition $\sigma'_S$ is still in $\Stab^{\dag}(S)$ and one has a birational map 
$$\modulis_{\sigma'_S}(\psi(0,0,1))\dashrightarrow\modulis_H(\psi(0,0,1))$$
coming from Theorem~\ref{K3-Stab_K3-large-volume-limit}. It will lead to an isomorphism $\modulis_H(\psi(0,0,1))\cong S$.

\begin{example}
The moduli space $\modulis_H(0,E,-2)$ is the relative compactified Jacobian $$\CJac_{-2}(S/|E|)\rightarrow |E|=\bp^1$$ 
which is in general not canonically isomorphic to $S$. In fact, one has $S\cong\modulis_H(0,E,-1)$.
\end{example}

\begin{example}
The moduli space $\modulis_H(1,L,-1)$ is canonically isomorphic to $S$. In fact, every element in it has the form $\shi_x(L)$ for some $x\in S$.
\end{example}

\subsection{The implications of the conjecture}
Now let us return to the EPW varieties. Assume that Conjecture \ref{Conjecture} is true, we are able to generalize Proposition~\ref{00_main-sextic-weak} and Corollary~\ref{00_main-dual-surface}.

\begin{proposition}\label{Conjecture_EPW-sextics-main}
Assume that Conjecture \ref{Conjecture} holds, the moduli space $\modulis_{\sigma_S}(v_2)$ is isomorphic to the double EPW sextic for any strongly smooth Gushel--Mukai surface $S$.
	
\begin{proof}
According to Proposition~\ref{EPW-sextics_generic} and \cite[Lemma 3.8]{KP18}, it remains to prove the statement when the special Gushel--Mukai threefold $X$ over $S$ admits a period dual $X'$ which is a special Gushel--Mukai threefold over $S'$. Using the conjecture, one finds an equivalence 
$$\Ku(X)\cong\Ku(X')$$
which sends $\kappa_1$ to $\kappa'_2$ and $\kappa_2$ to $\kappa'_1$ up to signs. It descends to an $\mathfrak{S}_2^\vee$-equivariant equivalence 
$$\DCoh(S)\cong\DCoh(S')$$
by \cite[Lemma 4.9]{BP23}, such that the stability condition $\sigma_S$ is sent to $\sigma_{S'}$ possibly up to a $\grp$-action and $v_2$ is sent to $v'_1$ up to a sign. Then one concludes by Proposition~\ref{EPW-sextics-dual_the-second-flop}.
\end{proof}
\end{proposition}

\begin{remark}
The parallel results for Gushel--Mukai fourfolds have been established in \cite{GLZ24,PPZ22} but only for the very general members. 
\end{remark}

Now let us prove Proposition~\ref{00_conjecture-generalization2}.

\begin{proposition}
Assume that Conjecture \ref{Conjecture} holds, then the double cover $\modulis_{\sigma_X}(\kappa_2)\rightarrow \modulis_{\sigma_S}(v_2)_{\mathfrak{S}_2^\vee}$ coincides with $\tilde{Y}^{\geq 2}_{A(X)}\rightarrow Y^{\geq 2}_{A(X)}$.
	
\begin{proof}
As in Proposition~\ref{Conjecture_EPW-sextics-main}, the case where $X$ has a special Gushel--Mukai threefold as its period dual can be proved using Proposition~\ref{EPW-surface-dual_special}. In general, one can take an ordinary Gushel--Mukai threefold $X'$ as the period dual according to Proposition~\ref{GM_threefold-period-dual}. Then the conjecture provides an equivalence
$$\Ku(X)\cong\Ku(X')$$
which sends $\sigma_X$ to a stability condition $\sigma_{X'}$ on $\Ku(X')$ and sends $\kappa_2$ to $\kappa'_1$. One can see 
$$\modulis_{\sigma_{X'}}(\kappa'_1)\cong\tilde{Y}^{\geq 2}_{A(X')^{\perp}}$$ 
as \cite[Theorem 7.3]{DK24} and \cite[Theorem 7.12]{JLLZ24} also apply to the ordinary case. Then one can conclude that the double cover $\modulis_{\sigma_X}(\kappa_2)\rightarrow \modulis_{\sigma_S}(v_2)_{\mathfrak{S}_2^\vee}$ is the one we want.
\end{proof}
\end{proposition}

As a result, one can easily conclude the following corollary using Proposition~\ref{En_image-forgetful-and-inflation} and the identifications of $\modulis_{\sigma_S}(v_1)$ and $\modulis_{\sigma_S}(v_2)$. It is exactly Corollary~\ref{00_conjecture-corollary}.

\begin{corollary}
Assume that Conjecture \ref{Conjecture} holds and let $X$ and $X'$ be two special Gushel--Mukai threefolds. Then one has a polarized isomorphism $\tilde{Y}^{\geq 1}_{A(X)}\cong \tilde{Y}^{\geq 1}_{A(X')}$ if and only if there is an equivalence $\Ku(X)\cong\Ku(X')$ sending $\kappa_1$ to $\kappa'_1$ and $\kappa_2$ to $\kappa'_2$ up to signs; one has a polarized isomorphism $\tilde{Y}^{\geq 1}_{A(X)}\cong \tilde{Y}^{\geq 1}_{A(X')^{\perp}}$ if and only if there is an equivalence $\Ku(X)\cong\Ku(X')$ sending $\kappa_1$ to $\kappa'_2$ and $\kappa_2$ to $\kappa'_1$ up to signs. 
\end{corollary}

\subsection{Description of the special period partners}
It is shown in \cite{JLZ26,Liu26} that the stable object $G_X$ representing the distinguished singular point $\tilde{p}_S$ in $\modulis_{\sigma_X}(\kappa_1)$ has the form 
$$\shu_X[2]\rightarrow G_X\rightarrow\shv_X^\vee[1]\rightarrow\shu_X[3]$$
where $\shu_X$ and $\shv_X$ are restrictions of tautological bundles on $\Gr(2,V_5)$ and $\Gr(3,V_5)\cong\Gr(2,V_5)$.

\begin{theorem}[{{\cite[Theorem 1.1]{JLZ26}}}]\label{BN-reconstruction}
Consider a special Gushel--Mukai threefold $X$ and the special stable object $G_X$ in $\modulis_{\sigma_X}(\kappa_1)$, then the Brill--Noether locus
$$\BN_X(G_X)=\{A\in\modulis_{\sigma_X}(\kappa_1+2\kappa_2)\,|\,\Hom(A,G_X[2])=\bc^{\oplus 3}\}\subset\modulis_{\sigma_X}(\kappa_1+2\kappa_2)$$
consists of the objects $K_x=\funct{P}\sho_x$ for any $x\in X$, where $\funct{P}$ is the left adjoint functor of the inclusion $\Ku(X)\subset\cate{D}^b(X)$. Moreover, one has $\BN_X(G_X)\cong X$.
	
\begin{proof}
One notices that the Kuznetsov component used in \cite{JLZ26} is different from $\cate{Ku}(X)$ but equivalent. The relation between the distinguished stable object $G_X$ and the characterizing object in \cite[Theorem 1.1]{JLZ26} is explained in \cite[Section 5]{JLLZ24}. The statement about the stable objects in the Brill--Noether locus follows from the argument for \cite[Theorem 1.1]{JLZ26}.
\end{proof}
\end{theorem}

The representatives for other singular points in the moduli space $\modulis_{\sigma_X}(\kappa_1)$ are determined in \cite[Corollary 6.11]{Liu26}. Each of them is given by a conic in the ramification locus of the covering involution on $X$ and has the form $\shi_C[1]$. Consequently, one can see

\begin{corollary}
Assume that Conjecture \ref{Conjecture} holds. Then a non-trivial special Gushel--Mukai threefold period partner of $X$ is isomorphic to the Brill--Noether locus
$$\BN_X(\shi_C[1])=\{A\in\modulis_{\sigma_X}(\kappa_1+2\kappa_2)\,|\,\Hom(A,\shi_C[3])=\bc^{\oplus 3}\}$$
in the moduli space $\modulis_{\sigma_X}(\kappa_1+2\kappa_2)$, where $C$ is a conic in the ramification locus of the branched double cover $X\rightarrow M_S$. This period partner is denoted by $X_C$.

\begin{proof}
Let $X'$ be such a period partner, then one can find an equivalence
$$\Ku(X)\cong\Ku(X')$$
using Conjecture \ref{Conjecture}, such that $\kappa_i$ is sent to $\kappa'_i$ for $i=1,2$. It follows that
$$\modulis_{\sigma_X}(\kappa_1)\cong\modulis_{\sigma_{X'}}(\kappa'_1)\quad\textup{and}\quad \modulis_{\sigma_X}(\kappa_1+2\kappa_2)\cong\modulis_{\sigma_{X'}}(\kappa'_1+2\kappa'_2)\quad$$ 
according to \cite[Appendix]{JLLZ24}. In particular, the singular point in $\modulis_{\sigma_{X'}}(\kappa'_1)$ represented by $G_{X'}$ is sent to a singular point $p$ in $\modulis_{\sigma_X}(\kappa_1)$. Suppose that $p$ is represented by $G_X$, then one has an isomorphism $X\cong X'$ by Theorem \ref{BN-reconstruction}, violating our assumption. So $p$ must be represented by $\shi_C[1]$ and one concludes that $X'\cong\BN_{X'}(G_{X'})$ is isomorphic to $\BN_X(\shi_C[1])$.
\end{proof}
\end{corollary}

One can see that $\BN_X(\shi_C[1])$ does not intersect $\BN_X(G_X)$ in $\modulis_{\sigma_X}(\kappa_1+2\kappa_2)$ as
$$\dim\Hom(\funct{P}\sho_x,\shi_C[3])=\dim\Hom(\sho_x,\shi_C[3])=\dim\Hom(\shi_C,\sho_x)<3$$
while it is expected that $X_C$ is the elementary transformation of $X$ along the conic $C$.\footnote{The reader can find the construction of elementary transformations in \cite[Section 4.4]{IP}.}

\section*{Acknowledgments}
The first author would like to thank his advisors Laura Pertusi and Paolo Stellari for their constant support, and to thank Arend Bayer and Jieao Song for helpful conversations. The authors would like to thank Kieran O'Grady and Sasha Kuznetsov for answering their questions. 

The authors also thank the anonymous referees for beneficial comments.

\end{document}